\documentclass[11pt]{amsart}
\usepackage{amsmath,amssymb}
\usepackage{amsfonts,amscd,latexsym}
\usepackage{euscript}
\newcommand{\labell}[1] {\label{#1}}

\newcommand{\sss}{{\smallskip}}
\newcommand{\T}{{\mathbb T}}
\newcommand{\Nn}{{\mathcal N}}
\newcommand{\Ff}{{\mathcal F}}

\newcommand{\tom}{{\widetilde \om}}
\newcommand{\tx}{{\widetilde x}}
\newcommand{\Th}{{\widetilde h}} \newcommand{\tg}{{\widetilde g}}

\newcommand{\TG}{{\widetilde G}}
 
\newcommand{\TF}{{\widetilde F}}
\newcommand{\ham}{{\it Htop}}

\newcommand{\wF}{{\widehat F}}
\newcommand{\Hams}{{\rm Ham}^{s}}
\newcommand{\Map}{{\rm Map}} 
 \newcommand{\Po}{{{\mathcal P}_{\om}}}
\newcommand{\Poz}{{{\mathcal P}_{\om}^{\Z}}}
\newcommand{\Poq}{{{\mathcal P}_{\om}^{\Q}}}
 
\newcommand{\Ppp}{{{\mathcal
P}}}

\newcommand{\Tg}{{\widetilde g}} 
 
\newcommand{\EM}{{\rm EM}} 
\newcommand{\TM}{{\widetilde M}} \newcommand{\Ta}{{\widetilde a}}
\newcommand{\TB}{{\widetilde B}} \newcommand{\Bb}{{{\mathcal B}}}
\newcommand{\Oo}{{{\mathcal O}}} \newcommand{\OM}{{{\mathcal O}^M}}
\newcommand{\OMG}{{{\mathcal O}^M_{\Ga}}} 
\newcommand{\OML}{{{\mathcal O}_{\La}^M}} 

 \newcommand{\id}{{\rm id}}
\newcommand{\MS}{{\medskip}} \newcommand{\Tor}{{\rm Tor}}

\newcommand{\C}{{\mathbb C}}
\newtheorem{theorem}{Theorem}[section]
\newtheorem{thm}[theorem]{Theorem}

\newtheorem{cor}[theorem]{Corollary}
\newtheorem{lemma}[theorem]{Lemma}

\newtheorem{prop}[theorem]{Proposition}
\newtheorem{defn}[theorem]{Definition}
\newtheorem{example}[theorem]{Example}
\newtheorem{rmk}[theorem]{Remark}
\newtheorem{question}[theorem]{Question}
\numberwithin{figure}{section}
\numberwithin{equation}{section}
\numberwithin{table}{section}
\newcommand{\p}{{\partial}}
\newcommand{\al}{{\alpha}}

\newcommand{\be}{{\beta}}

\newcommand{\Om}{{\Omega}}
\newcommand{\om}{{\omega}}
\newcommand{\eps}{{\varepsilon}}
\newcommand{\de}{{\delta}}

\newcommand{\ga}{{\gamma}}
\newcommand{\Ga}{{\Gamma}}
\newcommand{\io}{{\iota}}
\newcommand{\ka}{{\kappa}}
\newcommand{\la}{{\lambda}}
\newcommand{\La}{{\Lambda}}
\newcommand{\si}{{\sigma}}
\newcommand{\Si}{{\Sigma}}

\newcommand{\Kk}{{\mathcal K}}
\newcommand{\Ll}{{\mathcal L}}
\newcommand{\Aa}{{\mathcal A}}

\newcommand{\ov}{\overline}
\newcommand{\QED}{\hfill$\Box$\medskip}
\renewcommand{\Tilde}{\widetilde}
\newcommand{\inv}{{\rm inv}}

\newcommand{\Gg}{{\mathcal G}}

\newcommand{\Mo}{(M,\omega)}
\newcommand{\Ma}{(M,a)}
\newcommand\Hom{\operatorname{Hom}}

\newcommand\Flux{\operatorname{Flux}}
\newcommand\Symp{\operatorname{Symp}}
\newcommand\Ham{\operatorname{Ham}}

\newcommand\im{\operatorname{Im}}

\newcommand\Diff{\operatorname{Diff}}

\newcommand{\Hh}{{\mathcal H}}
\newcommand{\Q}{{\mathbb Q}}

\newcommand{\SO}{{\rm SO}}
\newcommand{\R}{{\mathbb R}}
\newcommand{\Z}{{\mathbb Z}}

\newcommand{\Tf}{{\widetilde f}}
\newcommand{\NI}{{\noindent}}
\begin{document}
    
\title{Enlarging the Hamiltonian group}
 
\author{Dusa McDuff}
\address{Department of Mathematics,
 Stony Brook University, Stony Brook, 
NY 11794-3651, USA}
\email{dusa@math.sunysb.edu}
\urladdr{http://www.math.sunysb.edu/~{}dusa}
\thanks{Partly supported by the NSF grant DMS 0305939.}
\date{6 January 2005, revised 1 May 2005 and 8 November 2005}
    \begin{abstract}  This paper investigates ways to enlarge the 
    Hamiltonian subgroup $\Ham$ of the symplectomorphism group $\Symp$ of a symplectic manifold $(M, \omega)$  to a group 
    that both intersects every connected component of $\Symp$ and
     characterizes symplectic bundles with fiber $M$ and closed 
    connection form.  
    As a consequence, it is shown that bundles with 
closed connection form are stable under appropriate small 
perturbations of the symplectic form.  
Further, the manifold $(M,\omega)$ has the property that every symplectic $M$-bundle has a closed connection form if and only if the flux group vanishes and the 
    flux homomorphism extends to a crossed homomorphism defined on the whole group $\Symp$.  The latter condition is equivalent to saying that a connected component of the commutator subgroup $[\Symp,\Symp]$ intersects the identity component of $\Symp$ only if it also intersects $\Ham$.      
    It is not yet clear when this condition is satisfied.
      We show that if the
    symplectic form vanishes on $2$-tori the flux homomorphism extends
    to the subgroup of $\Symp$ acting trivially on $\pi_1(M)$. 
    We also give an explicit formula for the Kotschick--Morita extension of Flux in the monotone case.  The results in this paper
    belong to the realm of soft symplectic topology, but raise some
    questions that may need hard methods to answer.
	\end{abstract}

\maketitle

   \tableofcontents

%%%%%%%%%%%%%%%%%%%%%%%%%%%%%%%%%%%%%%%%%%%%%%%%%%%%%%%%%%%%%%%%%%%%%%%%

\section{Introduction}\labell{sec:intro}
%%%%%%%%%%%%%%%%%%%%%%%%%%%%%%%%%%%%%%%%%%%%%%%%%%%%%%%%%%%%%%%%%%%%%%%%

%%%%%%%%%%%%%%%%%%%%%%%%%%%%%%%%%%%%%%%%%%%%%%%%%%%%%%%%%%%%%%%%%%%%%%%%
\subsection{Statement of the problem}\labell{ss:state}
%%%%%%%%%%%%%%%%%%%%%%%%%%%%%%%%%%%%%%%%%%%%%%%%%%%%%%%%%%%%%%%%%%%%%%%%

This paper studies what might be called the gross 
algebro-topological structure of the symplectomorphism group $\Symp(M,\om)$ of a closed connected symplectic manifold $(M,\om)$.  It is very well known that the identity component $\Symp_0(M,\om)$ of this group supports a continuous homomorphism called Flux that takes values in
the quotient $H^1(M;\R)/\Ga$ of the first de Rham cohomology group with the flux group $\Ga$. (Precise definitions are given in 
\S\ref{ss:extflux} below. By recent work of Ono~\cite{Ono} the group $\Ga$,
which is the image of $\pi_1(\Symp)$ under the flux, is
 a discrete subgroup of $H^1(M;\R)$.) 
 
The kernel of Flux is the Hamiltonian group
$\Ham(M,\om)$, whose elements are the time-$1$ maps of time-dependent 
Hamiltonian flows. Since $\Ham(M,\om)$ is perfect (see Banyaga~\cite{Ban}), it supports no 
nontrivial homomorphisms. Hence the only proper normal subgroups of 
$\Symp_0$ lie between $\Ham$ and $\Symp$. One of the motivations
for this paper is to understand the extent to which this structure in the identity component of $\Symp$ extends to the full group. 
What interesting normal (or approximately so) subgroups does $\Symp$ contain?  We shall tacitly assume that $H^1(M;\R)\ne 0$, since otherwise the
 whole discussion becomes trivial.

One can look at this question either algebraically or topologically. 
From an algebraic point of view the most important feature here is the flux homomorphism. 
Kotschick and Morita pointed out in~\cite{KM} that although Flux  does not in general extend to a homomorphism defined on $\Symp$, it does extend to a {\it crossed homomorphism} at least when $(M,\om)$
is monotone, for example a Riemann surface of genus $>1$.   
Therefore one interesting question is:
\begin{quote}{\it
does Flux always extend to a crossed homomorphism from $\Symp$
to $H^1(M;\R)/\Ga$?}
\end{quote} 
%whether such an extension always exists.  
So far, this is unresolved.
However, Proposition~\ref{prop:HHam1} gives equivalent geometric conditions for such an extension to exist.

 Kotschick and Morita's proof
 that Flux does extend in the monotone case is indirect. 
 We explain it at the beginning of \S\ref{sec:lift}.  Theorem~\ref{thm:TF} below gives 
an explicit construction for the extension,  which is based on looking at lifts of the action of the elements of $\Symp$ to an appropriate line bundle over the universal cover of $M$.
Much of our work was motivated by the search for similar explicit formulas.
For example we show in \S\ref{sec:Fs} that the composite of
$\Flux:\Symp_0\to H^1(M;\R)/\Ga$ with a suitable projection 
$H^1(M;\R)/\Ga\to H^1(M;\R)/\Po)$ always extends.

From a topological/geometric point of view, 
it turns out that these questions are very closely connected to
properties of  locally trivial symplectic $M$-bundles
$
(M,\om)\to P\to B.
$
A basic problem here is:
\begin{quote}{\it 
 to understand when 
the fiberwise symplectic class $a: = 
[\om]$  extends to a class $\Tilde a\in H^*(P;\R).$ }
\end{quote}
For example,
if $B$ is simply connected, $\Tilde a$ exists if and only if the 
structural group of the bundle reduces to $\Ham$.
 One reason to be interested in such extensions is that 
the class $\Tilde a$ may be used to define
characteristic classes that carry interesting information:
see K\c edra--McDuff~\cite{KMc} and Kotschick--Morita~\cite{KM2}. 
Another is that by work of Lalonde, McDuff and Polterovich~\cite{LMP,LM}
 bundles with structural group $\Ham$ have interesting cohomological properties. 
 For instance, 
for many bases $B$ the rational cohomology of $P$ splits as the tensor product of the cohomology of $B$ with that of $M$.  

One would like to know to what extent such results remain true for more general symplectic bundles. More precisely, for which subgroups $\Hh$ of $\Symp$ could one hope to extend the above results  to bundles with structural group $\Hh$?
Thus another motivating problem is: 
\begin{quote}
{\it to define a subgroup $\Hh$ of $\Symp$ with the property that 
the structural group of $P\to B$ reduces to $\Hh$ exactly when $a$ extends to $\Tilde a$.}
\end{quote}
 Such a group $\Hh$ can be thought of as an enlargement of $\Ham$ ---
 whence the title of this paper.  
We shall see that suitable groups $\Hh$ always exist; in \S2 we 
construct such $\Hh$ as kernels of extensions of homomorphisms related to the flux.  Proposition~\ref{prop:stab} shows that the stability properties of the Hamiltonian group that were proved in~\cite{LMP} extend to $\Hh$.  

The next question is to understand the obstruction to the reduction
of the structural group to $\Hh$; equivalently, when does the classifying map $B\to B\Symp$ of the bundle $P\to B$ lift to $B\Hh$?
If the Flux group $\Ga$ vanishes, then, as 
was pointed out in K\c edra--Kotschick--Morita~\cite{KKM}, Kotschick and Morita's argument may be used 
to show that this obstruction vanishes if and only if the
Flux homomorphism extends.  Theorem~\ref{thm:obs} gives the general result, which holds even when $\Ga\ne 0$.  

In the rest of this introduction we discuss each of the  main questions above in more detail.  In \S\ref{ss:res} we explain some examples and
also (in Remark~\ref{rmk:hsz}) some of 
the special features that arise when $[\om]$ is an integral class.
\MS

 \NI
{\bf Note.\,\,} This paper was first submitted before I saw 
K\c edra--Kotschick--Morita~\cite{KKM}, although most of their paper was completed earlier than mine. In this revision I have added a few remarks
to clarify the relation of their work to mine.  I have also reworked some arguments using insights from Gal--K\c edra~\cite{GalK}.  \MS

\NI
{\bf Acknowledgements.\,\,}  The author thanks Swiatoslav Gal, 
Jarek K\c edra, Jack Milnor, Leonid 
Polterovich and Zhigang Han for useful discussions.  
She also thanks Gal and K\c edra 
for showing her early drafts of their paper~\cite{GalK} and making
various helpful comments about earlier versions of this paper.
 In particular,
they helped streamline the proof of Lemma~\ref{le:obst} and 
correct some details in  Proposition~\ref{prop:main}.
Finally she wishes to thank the referees for  
helping to improve the presentation of the results discussed here.

%%%%%%%%%%%%%%%%%%%%%%%%%%%%%%%%%%%%%%%%%%%%%%%%%%%%%%%%%%%%%%%%%%%%%%%%
\subsection{Extending the fiberwise symplectic form}\labell{ss:groupext}
%%%%%%%%%%%%%%%%%%%%%%%%%%%%%%%%%%%%%%%%%%%%%%%%%%%%%%%%%%%%%%%%%%%%%%%%

Let $(M,\om)$ be a closed connected symplectic manifold, and consider 
a locally trivial symplectic $M$-bundle
$$
(M,\om)\to P\to B
$$
over a connected base $B$.
We want to understand when the fiberwise symplectic class $a: = 
[\om]$  extends to a class $\Tilde a\in H^*(P;\R).$
By Thurston's construction, this is equivalent to saying that 
%%D
the family of 
fiberwise symplectic forms $\om_b, b\in B,$ 
has a closed extension $\Om$ to $P$.
 (Here we assume without loss of generality
 that $P\to B$ is smooth.)  For short, we will often call the family
  $\om_b$ simply the fiberwise symplectic form.  
  
  This topic was first studied by
 Gotay, Lashof, Sniatycki and Weinstein in~\cite{GLSW} where they 
 showed that each
 extension $\Om$ of the fiberwise symplectic form
 gives rise to an Ehresmann connection on the bundle
 $P\to B$ whose horizontal spaces are the 
 $\Om$-orthogonals to the fibers. This connection has 
symplectic 
 holonomy iff  the  restriction of $\Om$ over the preimages of
arcs in the base is 
closed, and it has Hamiltonian holonomy round all contractible loops 
iff
$\Om$ is closed.
Therefore we shall call closed extensions of the fiberwise form
 {\bf closed connection forms}.
%  One reason for our interest in such extensions is that 
%the class $\Tilde a = [\Om]$ may be used to define
%characteristic classes that carry interesting information:
%see K\c edra--McDuff~\cite{KMc} and Kotschick--Morita~\cite{KM2}.  

If $H^1(M;\R) = 0$ then the Guillemin--Lerman--Sternberg (GLS) 
construction 
provides a closed connection form on every symplectic bundle.  
In the general case,
we are looking for a group homomorphism $\Hh\to \Symp(M,\om)$
such that an $M$-bundle $P\to B$ over a connected finite simplicial complex
has a closed connection form iff its 
classifying map $\phi:B\to B\Symp$ lifts to $B\Hh$.
For short we shall say that such a homomorphism $\Hh\to \Symp(M,\om)$
(or simply the group $\Hh$) has  the {\bf extension property}.
In particular, the group $\Symp$ itself has the extension property 
iff every symplectic $M$-bundle has a closed connection form.

The GLS construction also implies that 
a symplectic bundle $P\to B$  over a simply connected base $B$ has a 
closed 
connection form  if and only if its structural group  can be 
reduced to the Hamiltonian group $\Ham\Mo$.  
  These fibrations  are classified by maps $\phi:B\to B\Symp_0$ 
(where 
$\Symp_0$ denotes the identity component of the group $\Symp$), and
in this restricted case we may take $\Hh\to \Symp_0(M,\om)$ to be the 
inclusion 
$\Ham\hookrightarrow\Symp_0$.   Hence the desired group $\Hh$ should be 
understood as a generalization 
of the Hamiltonian group.

There are several natural candidates for $\Hh$.
Perhaps the most elegant approach is due to Seidel~\cite{Sei}, who 
considers 
a second topology on the symplectomorphism group 
called the
{\bf Hamiltonian topology} with basis consisting of the sets 
$gU$, for $g\in \Symp$ and $U$ open in $\Ham$.  We write 
$\Symp^{\ham}$
for the symplectomorphism group in this topology, reserving $\Symp$ 
to denote the same group but with its usual $C^{\infty}$-topology.  
The inclusion
$$
\Symp^{\ham}\to \Symp
$$ 
is obviously continuous, but is not a homeomorphism when 
$H^1(M;\R)\ne 
0$.  In particular, the identity component of $\Symp^{\ham}$ is the 
Hamiltonian group, not $\Symp_0$.   The following result is implicit 
in~\cite{Sei}, and holds by an easy application of the GLS
construction: see \S\ref{sec:Hams}.

\begin{prop}\labell{prop:sympham} The inclusion $\Symp^{\ham}\to 
\Symp$ has the extension property, i.e.
a symplectic $M$-bundle has a closed connection form iff its 
classifying 
map lifts to $B\Symp^{\ham}$.
\end{prop}

The group $\Symp^{\ham}$ is natural but very large.  For example its 
intersection with $\Symp_0$ has uncountably many components when 
$H^1(M;\R)\ne 0$.  We define in \S\ref{ss:extflux} below a closed 
subgroup 
$$
\Ham^s
$$
of $\Symp^{\ham}$
that still has the extension property, but has the homotopy
type of a countable CW complex.\footnote
{
A proof that $\Symp$ and $\Ham$ have the homotopy type of a 
 countable CW complex is sketched in 
McDuff--Salamon~\cite[9.5.6]{MS2}.}
Another advantage of this group is that it has an algebraic (rather 
than topological)
relation to $\Symp$, which makes it easier to understand the
homotopy fiber of the induced map $B\Ham^s\to B\Symp$.

The group $\Ham^s$ is a union of connected components of
$\Symp^{\ham}$.  It intersects every component of $\Symp$,
and when $H^1(M;\R)\ne 0$ intersects $\Symp_0$ in a countably
infinite number of components.  
Hence in general this subgroup  is not 
closed in $\Symp$.  However, as we point out in  
Remark~\ref{rmk:smallH}, if $\Ga\ne 0$ no closed subgroup of $\Symp$ can have the extension property.  

One reason for this is that $\Ham$ (which is closed in $\Symp$ by Ono~\cite{Ono}) does not
quite classify the set of $\Symp_0$-bundles for which $a$ extends. To explain this, it is useful to introduce the following definitions.
We suppose as before that
$M\to P\to B$ is a smooth bundle with base equal to 
 a finite dimensional (possibly open) connected
 manifold with finite homotopy type and   fiber a closed symplectic manifold. Recall also that a symplectic bundle
 has classifying map $\phi:B\to B\Symp_0$ iff it is symplectically trivial over the $1$-skeleton $B_1$ of the base.

\begin{defn}  \labell{def:mext}
We shall say that a (possibly disconnected) 
subgroup $\Hh_0$ of $\Symp_0$ has the {\bf 
 restricted extension property} if the following condition holds:
a $\Symp_0$-bundle $M\to P\to B$  has a closed connection form iff its
classifying map $B\to B\Symp_0$
lifts  to $B\Hh_0$. 
Similarly a subgroup $\Hh$ of $\Symp$ has the {\bf 
 modified extension property} if the following condition holds:
  a $\Symp$-bundle $M\to P\to B$  has a closed connection form iff 
  the pullback of its classifying map $B\to B\Symp$ over some finite cover
   $\rho:\TB\to B$ 
lifts  to $B\Hh$. 
\end{defn}
 
   Thus to say that $\Hh_0$ has the {\bf 
 modified restricted extension property} means that a
 $\Symp_0$-bundle $M\to P\to B$ has a closed connection form
iff there is a homotopy commutative diagram
$$
\begin{array}{ccc}  \Tilde B & \to & B\Hh_0\\
   \rho \downarrow && \downarrow\\
   B&\stackrel{\phi}\to& B\Symp_0,
   \end{array}
   $$
where $\rho:\TB\to B$ is some finite  covering map and $\phi:B\to B\Symp_0$ classifies $P\to B$.

One aim of this paper is to understand  subgroups $\Hh$ of $\Symp$  that 
have the (possibly modified) extension property. Since every $M$-bundle 
 $P\to S^1$ has a closed connection form, any  such group $\Hh$
 must intersect  almost every component of $\Symp$. 
The following proposition is proved  in \S\ref{ss:obs}.  We write  $\im(\pi_0(\Hh))$ for the image of $\pi_0(\Hh)$ in $\pi_0(\Symp)$.

  \begin{prop}\labell{prop:modext}  Let $\Hh$ be a subgroup of $\Symp$
with identity component equal to $\Ham$.
 Then $\Hh$ has the modified 
 extension property iff every finitely generated subgroup of
 $\pi_0(\Symp)$ has finite image in the coset space
$\pi_0(\Symp)/\im(\pi_0(\Hh))$.
 \end{prop}

 Lalonde--McDuff~\cite{LM} and early reprints of McDuff--Salamon~[Thm~6.36]\cite{MS} 
claim that the Hamiltonian group $\Ham\Mo$ has the  restricted 
extension property.
 But this is false: there are $\Symp_0$-bundles $P\to B$ that have a closed connection form but yet 
 only acquire a   Hamiltonian structure when pulled back over some 
 finite covering $\TB\to B$. (See later reprints of 
\cite{MS} and the erratum to~\cite{LM}.)  The  next 
proposition is proved in~\S\ref{ss:obs}.

\begin{prop}\labell{prop:ham}  The Hamiltonian group $\Ham\Mo$ has the 
 modified restricted  extension property.  It has the 
 restricted extension property iff $\Ga=0$.
 \end{prop}

The previous results prompt the 
 following question.

\begin{question}\labell{qu:2}  When does $\Symp$ have a subgroup $\Hh$
with  the modified extension 
property and such that 
$\Hh\cap \Symp_0=\Ham$?
\end{question}

As is shown by Proposition~\ref{prop:Flux} below, this is very closely related
to questions about extending
 the Flux homomorphism.  
%
%Since the group $\Ham$ is 
% closed in $\Symp$\footnote{This is 
%equivalent to the discreteness of the flux subgroup $\Ga$, a result 
%recently proved by Ono~\cite{Ono}.} one might hope to find a closed subgroup
%of $\Symp$ with the extension property.  But   
% if $\Ga\ne 0$ the Hamiltonian group 
%itself does {\it not}  have the extension property 
%appropriate to subgroups of $\Symp_0$: cf. 
%Proposition~\ref{prop:ham}.  Similar arguments show that when $\Ga\ne 0$ 
%no closed subgroup of $\Symp$ has the extension property, though there sometimes are 
% closed subgroups with the modified extension property of Definition~\ref{def:mext}. If these exist then
% one can define a smaller  subgroup than $\Hams$ with
%the extension property: see Remark~\ref{rmk:smallH}.
% 

The next problem is to understand 
the obstruction to the existence of a closed connection form.
The following lemma is proved in
 K\c edra--McDuff~\cite{KMc}.

\begin{lemma}\labell{le:2skel}  A symplectic $M$-bundle $\pi:P\to B$ 
has a closed connection 
form iff the restriction of $\pi$ over the $2$-skeleton of $B$ has 
such a form.
\end{lemma}

This  is mildly surprising:  in order for the 
fiberwise symplectic class $[\om]$ to extend to $H^2(P)$ it must lie 
in 
the kernel of the Leray--Serre  differential $d_3$ 
 as well as in ${\rm ker\,}d_2$, and in principle $d_3$ depends on 
the $3$-skeleton of $B$.  However, Lemma~\ref{le:2skel} 
is a very general result that 
is valid in the cohomologically symplectic case, i.e. for 
pairs $\Ma$ where $M$ is a closed oriented $2n$-manifold and $a\in 
H^2(M)$ has $a^n>0$.  To prove it, observe that if $a$ survives to 
$E_{3}^{0,2}$ then $d_3(a^{n+1}) = (n+1)d_3(a)\otimes a^n$ must vanish
since $a^{n+1}=0$.  But because tensoring with $a^n$ gives an 
isomorphism $E_3^{3,0}\to E_3^{3,2n}$ and the cohomology groups have 
coefficients $\R$,  this is possible only if 
$d_3a=0$.

Although the obstruction lies in such low dimensions, it is still not 
fully understood.  It divides into two parts, one that depends on the \lq\lq symplectic mapping class group"
$\pi_0(\Symp)$ and the other on the flux subgroup $\Ga$.
% (whose definition is 
%recalled in \S\ref{ss:extflux} below).
This is shown by the next proposition that
formulates  necessary and
sufficient conditions for the obstruction to vanish.
% in the case
%of bundles where $\pi_1(B)$ acts trivially on 
%$H^1(M;\R)$.  
% but this condition itself is not well understood: see
%Remark~\ref{rmk:II}.
%It is given in terms of the flux group $\Ga$ (whose definition is 
%recalled in \S\ref{ss:extflux} below) and the  
%We write $\Symp_H$ for the group of all symplectomorphisms that act
%trivially on $H^1(M;\R)$, and  $[G,G]$  for the commutator
%subgroup of $G$.

\begin{prop}\labell{prop:noact} The following conditions are
    equivalent:
    
    \NI
    {\rm (i)}  Every symplectic $M$-bundle $P\to B$ 
has a closed connection form.
\smallskip

\NI
{\rm (ii)}  $\Ga = 0$ and every connected component of the commutator 
subgroup $[\Symp,\Symp]$ that intersects $\Symp_0$ also intersects $\Ham$.
\end{prop}

The second condition in (ii) is not yet well understood. 
We show below that it  is equivalent to
the existence of a suitable extension of the flux homomorphism;
see  Proposition~\ref{prop:HHam1} and Remark~\ref{rmk:II}.
However, we can prove that Flux extends only under very restrictive circumstances, for example if  $[\om]$ vanishes on $2$-tori and  $\pi_0(\Symp)$ acts
 on $\pi_1(M)$ by inner automorphisms; see  Proposition~\ref{prop:symphyp}.

%%%%%%%%%%%%%%%%%%%%%%%%%%%%%%%%%%%%%%%%%%%%%%%%%%%%%%%%%%%%%%%%%%%%%%%%
 
\subsection{Extending the flux homomorphism}\labell{ss:extflux}
%%%%%%%%%%%%%%%%%%%%%%%%%%%%%%%%%%%%%%%%%%%%%%%%%%%%%%%%%%%%%%%%%%%%%%%%

The flux is initially defined as a homomorphism from the
universal cover $\Tilde{\Symp}\,\!_0$ of the identity component 
of the symplectomorphism group to the group
$H^1(M;\R)$.  For each element $\{g_t\}\in \Tilde{\Symp}\,\!_0$ the value of 
the class $\Tilde\Flux(\{g_t\})\in H^1(M;\R)$ on the $1$-cycle $\ga$ in $M$ is given by integrating $\om$ over the $2$-chain $(s,t)\mapsto g_t(\ga(s))$.
If we define the Flux group $\Ga$ to be the 
image of $\pi_1(\Symp)\subset \Tilde{\Symp}\,\!_0$  under  $\Tilde\Flux$, then $\Tilde\Flux$ descends to a homomorphism 
$$
\Flux: \Symp_0\Mo\to H^1(M;\R)/\Ga
$$
that we shall call the {\bf flux homomorphism}.
Its kernel is the Hamiltonian group $\Ham$.

The problem of extending the Flux homomorphism 
to the whole group $\Symp$ arose  (with rather different motivation) 
in the work of
Kotschick--Morita~\cite{KM} in the case  
when $M$ is a  Riemann surface of genus $g>1$ or, more generally a 
monotone manifold, i.e. a manifold in which 
the symplectic class $[\om]$ is a multiple 
of the first Chern class.
 They showed that in this case Flux extends to a 
{\bf crossed homomorphism}
$$
F_{KM}: \Symp(M,\om)\to H^1(M;\R),
$$
that is, a map $F: = F_{KM}$ that
instead of being a homomorphism, satisfies the identity
\begin{equation}\labell{eq:crh}
F(gh) = F(h) + h^* F(g),
\end{equation}
where $h^*$ denotes the action of $h$ on $H^1(M;\R)$
via pullback.\footnote
{
This is the natural identity to use for a crossed homomorphism $G\to 
\Aa$  when the group $G$ acts 
contravariantly on the coefficients $\Aa$. Note also that
when $\Symp$ acts nontrivially on $H^1(M;\R)$ it is not
possible to extend Flux to a group homomorphism: see 
Remark~\ref{rmk:comm}.}

In general, one should look for an extension of Flux
with values in $H^1(M;\R)/\Ga$.
So far, it is unknown whether an extension must always exist: see
Proposition~\ref{prop:HHam1}.
However, the following result shows that this question is very 
closely related to 
our earlier considerations.

\begin{prop} \labell{prop:Flux}  {\rm (i)}\, If 
$$
\Tilde F:\Symp(M,\om)\to H^1(M;\R)/\Ga
$$
 is a continuous crossed homomorphism 
that extends $\Flux$,
its kernel $\Hh$ intersects every component of $\Symp$ and
has the modified extension property. Moreover $\Hh$ has the extension 
property iff $\Ga=0$.
\sss

\NI
{\rm (ii)}\,  Conversely, let  $\Hh$ be a subgroup of $\Symp$
that intersects $\Symp_0$ in $\Ham$
and denote by
 $\Symp_{\Hh}$ the union of the components of $\Symp$ that intersect 
 $\Hh$.  
Then there is a crossed homomorphism
$\TF:\Symp_{\Hh}\to  H^1(M;\R)/\Ga$ that extends $\Flux$. 
\end{prop}

\begin{proof}
   Given $\TF$, let  $\Hh: = \Hh_{\TF}$ 
   be the kernel of $\TF$.  Then $\Hh\cap\Symp_0 = \Ham$.  Further 
   given any 
   $g\in \Symp$ choose $h\in \Symp_0$ such that $\Flux h = -\TF(g)\in 
   H^1(M;\R)/\Ga$.  Then $g$ is isotopic to the element 
   $gh\in \Hh$.  Hence $\Hh$ has the modifed extension property by
   Proposition~\ref{prop:modext}.  If $\Ga=0$ then the inclusion 
   $B\Hh\to B\Symp$ is a homotopy equivalence and every bundle has 
   both a closed connection form and an $\Hh$-structure.  If $\Ga\ne 
   0$ one can construct bundles that have a closed connection form 
   but no $\Hh$-structure as in the proof of 
Proposition~\ref{prop:ham}.

 To prove (ii) we define  $\TF$ on $\Symp_{\Hh}$ as follows:
given $g\in \Symp_{\Hh}$  denote by 
$\si_g$ any element in $\Hh$ that is isotopic to $g$ and set
$
\TF(g): = \Flux(\si_g^{-1} g).
$
This is independent of the choice of $\si_g$.  Further
$
(\si_g\si_h)^{-1}\si_{gh} \in \Hh\cap \Symp_0 = \Ham.
$
Hence 
  $$
  \TF(gh) = \Flux\bigl((\si_{gh})^{-1} gh\bigr) = 
  \Flux\bigl(\si_h^{-1}(\si_g^{-1}g)\si_h\bigr) \Flux\bigl(\si_h^{-1}h) = h^*\TF(g) + \TF(h).
  $$
Thus $\TF$ 
satisfies~(\ref{eq:crh}) and so is a crossed homomorphism.
\end{proof}

\begin{rmk}\labell{rmk:smallH}\rm  If Flux extends to $\TF$ but $\Ga\ne 0$ then by part (i) of the above proposition the 
kernel of $\TF$ does not have the extension property. On the other hand, the kernel $\Hh_\Q$ of the composite map 
$$
\TF:\Symp\;\to\; H^1(M;\R)/\Ga\;\to\; H^1(M;\R)/(\Ga\otimes \Q)
$$
does have the extension property by Corollary~\ref{cor:main2}.
It follows from the proof of Proposition~\ref{prop:ham} 
that this is the smallest group with this property.  Note that it has countably many components in $\Symp_0$.
\end{rmk}

Although Flux may not always  have an extension with values 
in $H^1(M;\R)/\Ga$,
%%D(For further discussion of this issue, see 
%Remarks~\ref{rmk:II} and \ref{rmk:hsz}.)
 its composite with 
projection onto a suitable  quotient group $H^1(M;\R)/\La$ always can 
be extended.  Below we shall define a
continuous crossed homomorphism
\begin{equation}\labell{eq:aaa}
\widehat F_s: \Symp(M,\om)\to H^1(M;\R)/H^1(M;\Po) =:\Aa,
\end{equation}
where $\Po: = \Poq$ is the rational period group of $\om$ (i.e. the values 
taken by $[\om]$ on the rational $2$-cycles).   The map
$\widehat F_s$ depends on the choice of a splitting $s$ of a certain exact sequence. (See 
the definitions in \S\ref{sec:Fs}.  The \lq\lq topology" on $\Aa$ is explained in Remark~\ref{rmk:top}.)  However
its restriction  
to the identity component $\Symp_0$ is independent of this choice and
equals the composite
$$
\Symp_0\stackrel {\Flux}\to H^1(M;\R)/\Ga \to H^1(M;\R)/H^1(M;\Po).
$$
%(For more details see 
%Proposition~\ref{prop:wf} and Remark~\ref{rmk:top}.)

Recall that if a  group $G$ acts continuously on an $R$-module  $\Aa$ 
(for suitable ground ring $R$) then
the continuous group cohomology\footnote
{
The group cohomology of a discrete group $G^\de$ is defined algebraically using a cochain complex that in older literature is called the 
Eilenberg--MacLane complex.  This cohomology is isomorphic to
the singular cohomology of the classifying space $BG^\de$. 
 If $G$ is a topological group, then one can also consider the cohomology groups defined by restricting to continuous cochains.
% 
%a (continuous) 
%group cohomology defined using the (continuous) Eilenberg--MacLane complex.  
Since the resulting cohomology groups are completely  different from the singular cohomology of 
$BG$, we will for the sake of clarity denote the group cohomology by $H^*_{\EM}$, adding a $c$ wherever appropriate to emphasize continuity.
}
$$
H_{c\EM}^1(G; \Aa)
$$
(defined using continuous Eilenberg--MacLane cochains) is the 
quotient of 
the module of all continuous crossed homomorphisms $G\to \Aa$
by the submodule consisting of the coboundaries $h\mapsto h\cdot\al - 
\al, \al\in \Aa$.  Therefore, 
$\wF_s$  defines an element
$$
[\wF_s]\in H_{c\EM}^1(\Symp; \Aa).
$$
Although there is no canonical choice for $\wF_s$ it turns out that 
the cohomology class
$[\wF]: = [\wF_s]$ is independent of the choice of $s$. 
 
We now define
$$
\Hams\Mo: = \ker \wF_s.
$$
These groups depend on the representative $\wF_s$ chosen for the 
class
$[\wF_s]$, but they are all  conjugate via elements of $\Symp_0$.
Moreover their 
intersection with the subgroup $\Symp_H$ of $\Symp$ that acts 
trivially on $H^1(M;\R)$ is independent of $s$. (See 
Lemma~\ref{le:conj}.)  This holds because  any
crossed homomorphism $\wF: \Symp\to \Aa$ restricts to a  
homomorphism on  $\Symp_H$ that depends only on the class represented 
by 
$\wF$ in  $H_{c\EM}^1(\Symp; \Aa)$.

    Because $\Hams$ is the kernel of a continuous crossed homomorphism it 
    follows from standard theory that one can use this 
    homomorphism to define a  class 
    $$
{\OM}\in H^2(B\Symp; H_\Q)
$$
that measures the obstruction to lifting a map $\phi:B\to B\Symp$ to 
    $B\Hams$. Here one should think of $H_\Q:=H^1(M;\Po)$ as
$\pi_1(\Aa)$ where  $\Aa: = H^1(M;\R)/H^1(M;\Po)$.
%
% as a 
%quasitopological group: see Remark~\ref{rmk:top}.  Thus 
%the local coefficient system $\pi_1(\Aa)$ on $B\Symp$ 
%has fibers isomorphic to the discrete 
%group $H_\Q: = H^1(M;\Po)$.  

The next result explains the role of this class $\OM$.
The statement in (i) below 
arose from some remarks in Gal--K\c edra~\cite{GalK}.

\begin{thm}\labell{thm:obs}
    {\rm (i)} The obstruction class
     $\OM\in H^2(B\Symp;H_\Q)$ equals the image  $d_2^\om([\om])$
     of $[\om]\in H^2(M;\Po)$ under the differential $d_2^\om$ in the 
     Leray--Serre spectral sequence for the cohomology
of the universal $M$-bundle over $B\Symp$ with coefficients $\Po$.    \smallskip    
    
    \NI
    {\rm (ii)}   There is a  crossed homomorphism $\Tilde 
    F: \Symp\to H^1(M;\R)/\Ga $ that extends Flux if and only if
    $\OM$ lies
    in the image of $ 
    H^2(B\Symp;\Ga)$ in $H^2(B\Symp;H_\Q)$.\sss
%    
%    \NI
%    {\bf (iii)}  If $\Ga=0$ and the hypothesis in (ii) holds, then  
%$\OM=0$. 
\end{thm}

\begin{cor}\labell{main}
 $\Hams$ has the extension property.
\end{cor}
\begin{proof}  By Lemma~\ref{le:2skel} there is a closed extension
of $\om$ iff $d_2([\om])=0,$ where $d_2$ denotes the differential in the spectral sequence for real cohomology.  But because $\Po$ is divisible, this vanishes iff $d_2^\om([\om])=0.$
\end{proof}

\begin{cor}\labell{cor:obs} The  following  conditions
are equivalent:
\sss

\NI
{\rm (i)}  Every symplectic $M$-bundle has a closed connection form. 
\sss

\NI
{\rm (ii)}
$\Ga=0$ and there is a crossed homomorphism $\TF:\Symp\to H^1(M;\R)$ extending Flux.
\end{cor}

For example, when $\Mo$ is  monotone, 
the first Chern class of the vertical tangent bundle of $P\to B$ 
provides an extension of $[\om]$.   Therefore 
the obstruction class $\OM$ must vanish.  This is consistent with 
the corollary 
since the Kotschick--Morita homomorphism 
$F_{KM}$ extends Flux. 

The next result clarifies the conditions under which  Flux can be extended.  In Kotschick--Morita~\cite[\S6.3]{KM}, the obstruction to 
the existence of an extension of Flux is described as
a certain class $\de([F])\in H^2_{\EM}\bigl(\pi_0(\Symp); H^1(M;\R)/\Ga\bigr)$.  Thus the next proposition can be interpreted as giving  geometric explanations of what it means for this class to vanish.

\begin{prop} \labell{prop:HHam1} The  following  conditions
are equivalent.
\sss

\NI
{\rm (i)} $\Flux$ extends to a crossed homomorphism
$
\Tilde F:\Symp(M,\om)\to H^1(M;\R)/\Ga.
$
\sss

\NI
{\rm (ii)}
    For every closed Riemann surface 
$\Si$ every representation of $\pi_1(\Si)$ in 
$\pi_0(\Symp)$ lifts to a representation into the group $\Symp/\Ham$.
\sss

\NI
{\rm (iii)}
 For every product of commutators $[u_1,u_2]\dots[u_{2p-1},u_{2p}]$, $u_i\in \Symp$, that lies in $\Symp_0$, there are elements $g_1,\dots,g_{2p}\in \Symp_0$ such that 
$$
[u_1g_1,u_2g_2]\dots[u_{2p-1}g_{2p-1},u_{2p}g_{2p}] \in \Ham.
$$

\NI
{\rm (iv)} For every  symplectic $M$-bundle $P\to \Si$ there is a 
    bundle $Q\to S^2$  such that the fiberwise connect sum $P\#Q\to 
    \Si\#S^2 = \Si$ has a closed connection form. 
    \end{prop}

\begin{rmk}\labell{rmk:II}\rm
     If we restrict to the subgroup
    $\Symp_H$  of $\Symp$ that acts trivially on $H^1(M;\R)$
then  (iii) is equivalent to saying that
\begin{equation}\labell{eq:commut5}
[\Symp_H, \Symp_H] \cap \Symp_0 = \Ham.
\end{equation}
   We show in Corollary~\ref{cor:ext} 
    that equation (\ref{eq:commut5}) holds
     when $[\om]$ vanishes on tori and  lies in the subring of $H^*(M)$ generated by $H^1$.  However,
it is so far unknown whether it always holds.  If not, then $\TF$
    cannot always exist.  On the other hand, there are indications that
    (\ref{eq:commut5}) might always hold.  It seems that a large part of $\pi_0(\Symp)$ can be
    generated by Dehn twists about Lagrangian spheres: cf. 
    Seidel~\cite[1.7]{Sei}.  In dimensions $>2$ these are well defined up
    to Hamiltonian isotopy and act trivially on $H^1(M)$, and so one might
    be able to take $\Hh\cap \Symp_H$ to be the group generated by Dehn
    twists.  In any case, it does not seem that the methods used in this
    paper are sufficiently deep to resolve this question.
\end{rmk}

%%%%%%%%%%%%%%%%%%%%%%%%%%%%%%%%%%%%%%%%%%%%%%%%%%%%%%%%%%%%%%%%%%%%%%%%
 
\subsection{Further results and remarks}\labell{ss:res}
%%%%%%%%%%%%%%%%%%%%%%%%%%%%%%%%%%%%%%%%%%%%%%%%%%%%%%%%%%%%%%%%%%%%%%%%

After discussing stability, we describe a few cases 
where it is possible to extend the Flux homomorphism.  We end by 
discussing the integral case, and the question of uniqueness.\MS

\NI {\bf Stability under perturbations of $\om$.}
It was shown in Lalonde--McDuff~\cite{LM} that Hamiltonian
bundles are stable under small perturbations of $\om$. One cannot 
expect general
symplectic bundles to be stable under arbitrary small perturbations of
$\om$ since $\pi_1(B)$ may act nontrivially on $H^2(M;\R)$.
Given a symplectic bundle $(M,\om)\to P\to B$ let us denote by 
$V_2(P)$  the subspace of $H_2(M;\Q)$ generated by 
the elements $g_*(C)- C$, where $C\in H_2(M;\Z)$ and $g$ is any symplectomorphism of $M$ that occurs as 
 the holonomy of a symplectic connection on $P\to B$ around some loop in $B$.
(Since $g_*(C)$ depends only on the smooth isotopy class of $g$, it 
does not matter which connection we use.)    The 
 subspace $H^2(M;\R)^{\inv(P)}$, consisting of
 classes $a\in H^2(M;\R)$ that are fixed by all such $g$,
 is the annihilator  of $V_2(P)$.  The 
most one can expect is that the existence of a symplectic structure 
on $P\to B$ is stable under perturbations of $[\om]$ in this subspace. 
For example, if
$\om$ is {\bf generic} in the sense that it gives an injective map
 $H_2(M;\Z)/{\rm Tor}\to \R$  then $V_2(P)$ is torsion and 
 $H^2(M;\R)^{\inv(P)} =H^2(M;\R)$.

 \begin{prop}\labell{prop:stab}  Let $(M,\om)\to P\to B$ be 
a symplectic $M$-bundle over a finite simplicial complex $B$.
Then there is a neighborhood $\Nn(\om)$ of $\om$ in
the space of all closed $2$-forms on $M$ that 
represent a class in $H^2(M;\R)^{\inv(P)}$
such that for all $\om'\in \Nn(\om)$:\smallskip

\NI
{\rm (i)} $P\to B$ has the structure of an $\om'$-symplectic bundle, 
and
\smallskip

\NI
{\rm (ii)}  if there is a closed extension of $\om$, then the same is 
true for $\om'$. 
\end{prop}

 Part (i) of this proposition follows by the arguments in 
\cite[Cor.~2.5]{LM}.  Part (ii) was also proved in~\cite{LM} in the 
case when $B$ is classified by a map into $B\Symp_0(M,\om)$.
The proof of the general case
is given at the end of \S\ref{sec:Hams}.  The next corollary 
is an immediate consequence of (ii).

\begin{cor}  If $P\to B$ has a $\Ham^s$-structure then the image of 
the restriction map $H^2(P;\R)\to H^2(M;\R)$ is the subspace 
$H^2(M;\R)^{\inv(P)}$ of $H^2(M;\R)$ that is invariant under 
the action of $\pi_1(B)$. \end{cor}

This result implies that the differentials $d_2^{2,0}$ 
and $d_3^{2,0}$ in the Leray--Serre 
spectral sequence for the real cohomology  of $P\to B$ vanish,
 and so is a partial generalization of the vanishing 
results in~\cite{LM}.

\begin{rmk}\labell{rmk:stab2}\rm
Proposition~\ref{prop:stab} is proved using the Moser homotopy argument
 and so works
only over compact pieces of $\Symp$.  This is enough to give stability for  bundles over  bases $B$ of finite homotopy type
 but is not enough to allow one to make any statements 
about properties that involve the full group $\Symp$.  Hence even if
$\Flux$ extends to $\TF:\Symp\to H^1(M;\R)/\Ga$ for the manifold $\Mo$,
it is not clear that it also extends for sufficiently close forms $\om'$ 
whose cohomology class is invariant under $\Symp(M,\om)$.  For one thing,
however close $\om'$ is, there may be new components of $\Symp(M,\om')$ containing elements that are 
far from those in $\Symp(M,\om)$. 
\end{rmk}
\MS

\NI
{\bf Manifolds with $\Ga=0$.}\,\,  One expects that for most manifolds
$\Ga = 0$.  Whether $\OM$ then vanishes
is still not clear.  We now discuss some special cases in which $\Ga
=0$ and Flux extends to
a crossed homomorphism defined either on the whole group $\Symp$
or on some large subgroup.

The first case is
when $\Mo$ is (strongly) monotone, i.e. the symplectic class $[\om]$
is a multiple of the first Chern class.  In this case $\OM=0$ since
one can always choose a closed connection form in the class of a multiple of the
vertical first Chern class.   Kotschick--Morita~\cite{KM} observed that
Flux always extends.  We shall give an explicit formula for $\TF$ in
Theorem~\ref{thm:TF}.  As 
noted in Remark~\ref{rmk:exten} (ii), the argument in 
fact applies whenever $[\om]$ is integral and $\Symp$ has the integral extension property, i.e. there is a complex line bundle over
the universal $M$-bundle $M_{\Symp}\to B\Symp$ whose first Chern class 
 restricts to $[\om]$ on the fiber.  

Another somewhat tractable case is when $\Mo$ is {\bf atoroidal}, i.e.
$\int_{\T^2} \psi^*\om = 0$ for all smooth 
maps $\psi:\T^2\to M$.  Note that $\Ga=0$ for such manifolds, because
for each loop $\{f_t\}$ in $\Symp_0$ the value of the
class $\Flux(\{f_t\})$ on the $1$-cycle $\ga$ is obtained by
integrating $\om$ over the torus $\cup_{t}f_t(\ga)$. 
In the next proposition, we denote by $\Symp_\pi$  the subgroup of $\Symp$ consisting of elements that are isotopic to a symplectomorphism that fixes the basepoint $x_0$ of $M$ and induces the identity map $\pi_1(M,x_0)\to \pi_1(M,x_0)$.\footnote
{
One can check that $g\in \Symp_\pi$ iff for any path $\ga$ in $M$ from $x_0$ to $g(x_0)$ the induced maps $\ga_*,g_*:\pi_1(M,x_0)\to \pi_1(M,gx_0)$ differ by an inner automorphism. Thus, loosely speaking, $\Symp_\pi$ consists of all symplectomorphisms that act trivially on $\pi_1(M)$.}

\begin{prop}\labell{prop:symphyp}  
 If $\Mo$ is atoroidal then $\Ga= 0$  and $\Flux$ extends to a 
homomorphism $\TF: \Symp_{\pi}\to H^1(M;\R)$.
\end{prop}

We shall see in \S\ref{ss:ator} that in the above situation $\TF$ can be extended to a crossed homomorphism defined on the whole of $\Symp$ but at the cost of enlarging the target group.  

Proposition~\ref{prop:symphyp}  gives a partial answer to K\c edra--Kotschick--Morita's question~\cite{KKM} of whether the usual flux homomorphism
$\Symp_0\to H^1(M;\R)$ extends to the full group $\Symp$ when $[\om]$
is a bounded class. This condition means that $[\om]$  may be represented by a
singular cocycle that is uniformly bounded on the set of all singular
$2$-simplices.\footnote { The (smooth) cocycle represented by
integrating $\om$ can never be bounded because bounded cocycles vanish
on cylinders as well as tori.} If $[\om]$ is bounded, then $\Mo$ is atoroidal 
since an arbitrary multiple of a toric class $C$ can be represented by the sum of
just two singular $2$-simplices.  Another interesting atoroidal case is that of
 symplectically hyperbolic manifolds.  There are various possible definitions here.  We shall use 
 Polterovich's definition from~\cite{Po} in which $\Mo$ is called
 {\bf symplectically hyperbolic} 
 if the pullback 
$\tom$ of $\om$ to the universal cover $\TM$ of $M$ 
has bounded primitive, i.e. $\tom = d\be$ for some $1$-form $\be$ 
that is bounded with respect to any metric on $\TM$ that is pulled 
back from $M$. For example, $\Mo$ might be a product of Riemann 
surfaces of genus $>1$ with a product symplectic form.
Because in the covering $\R^2\to \T^2$ the boundary of a square of 
side $N$ encloses $N^2$ fundamental domains, it is easy to check that
any $2$-form on $\T^2$ whose pullback 
to  $\R^2$ has bounded primitive must have zero integral over 
$\T^2$.  
Hence we find:

\begin{lemma}\labell{le:symphyp}  
Proposition~\ref{prop:symphyp} applies both when $[\om]$ is bounded
and when $(M,\om)$ is 
    symplectically hyperbolic.
    \end{lemma}
    
\S\ref{ss:ator} contains a few other similar results that are valid in
special cases, for example when $\om$ vanishes on $\pi_2(M)$. 

We end
the introduction with some general remarks.

\begin{rmk}\labell{rmk:hsz}\rm {\bf (The integral case.)}  There is
an analogous group $\Ham^{s\Z}$ which is the kernel of a crossed homomorphism
$\wF_s^\Z$ with values in $H^1(M;\R/\Poz)$ where $\Poz$ denotes the
set of values of $[\om]$ on the integral $2$-cycles $H_2(M;\Z)$.  In
many respects the properties of this group are similar to those of
$\Ham^{s}$.  However, there are some interesting differences.
If $\Tor$ denotes the torsion subgroup of $H_1(M;\Z)$, then the analog
of the group $\Aa$ occurring in equation~(\ref{eq:aaa}) is
$$
\Aa^\Z: = \;H^1(M;\R/\Poz)\;\cong \;\Hom(\Tor, \R/\Poz) \;\oplus\; H^1(M;\R)/H^1(M;\Poz).
$$
 Hence Theorem~\ref{thm:obs} (i) does not immediately generalize; the proof of Lemma~\ref{le:obst} shows that
the obstruction to the existence of a $\Ham^{s\Z}$-structure is twofold,
the first coming from the finite group $\pi_0(\Aa^\Z)$ (see Lemma~\ref{le:hsz0})
 and the second  an obstruction cocycle similar to $\OM$ coming from $\pi_1(\Aa^\Z)$ (see Lemma~\ref{le:OMobst}).
  Nevertheless,  
since every $\Hams$-bundle over a compact base $B$ has a finite cover
with a $\Ham^{s\Z}$-structure, the latter group has the modified extension
property.

The group $\Ham^{s\Z}$ is 
most interesting in the case when $\Poz=\Z$, i.e. when $[\om]$ is a 
primitive integral class.  In this situation one might expect 
$B\Ham^{s\Z}$ to classify bundles 
$P\to M$ that have a closed  and {\it integral}
connection form.  Even in the case when 
$\pi_0(\Symp)$ acts trivially on 
$H_1(M;\Z)$, this is not quite true.  Gal--K\c edra~\cite{GalK} point out
that there is a further torsion obstruction in $H^3(B;\Z)$ which measures whether the symplectic class $[\om]$ has an integral rather than rational extension: cf.
 Proposition~\ref{prop:main} below and Example~\ref{ex:moreinteg}.
% shows that there is a further torsion obstruction in $H^3(B;\Z)$
% this is the case when 
%$\pi_0(\Symp)$ acts trivially on 
%$H_1(M;\Z)$.  Gal--K\c edra~\cite{GalK} show that this remains true 
%when  $\pi_0(\Symp)$ acts trivially on 
%the torsion classes in $H_1(M;\Z)$.  (This hypothesis is equivalent 
%to saying that the action of $\Symp$ on the set of prequantum line 
%bundles is trivial.)

If 
$\pi_0(\Symp)$ acts nontrivially on $\Tor\subset
H_1(M;\Z)$, then the groups $\Ham^{s\Z}$ are not all conjugate. We show in Lemma~\ref{le:tau} that 
up to homotopy the choice of splitting $s$ is equivalent to 
the choice of an integral lift $\tau$ of
$[\om]$, i.e. of a prequantum line bundle $L_\tau$.  Moreover $s$ itself is determined by a unitary connection  $\al$ on $L$, and the group 
$\Ham^{s\Z}$ consists of all
symplectomorphisms $\phi$ that preserve the monodromy of $\al$, i.e. 
for all closed loops $\ga$ in $M$ the $\al$-monodromy  $m_\al(\ga)$ round $\ga$ 
equals that round $\phi(\ga)$.  Thus, $\Ham^{s\Z}$ is the same as the group $D_\ell$ 
considered by Kostant in~\cite{Ko}.  (It is also 
homotopy equivalent to the covering group of $\Symp$
 considered in Gal--K\c edra~\cite{GalK};  
see Proposition~\ref{prop:Z} below.) 
Hence one can think of the monodromy $m_\al$ of $\al$ as 
a \lq\lq Hamiltonian structure" on $M$, i.e. this function on the space $\Ll M$ of  closed loops in $M$ is the structure on $M$ that 
is preserved by the elements of $\Ham^{s\Z}$. 
\footnote
{
The function $m_\al: \Ll M\to \R/\Z$ is
characterized by the following two properties:
 (i)  
$m_\al(\be*\ga) = m_\al(\be) + m_\al(\ga)$, where $\be*\ga$ is the 
concatenation of two loops with the same base point;
and
(ii)  if $\ga$ is the boundary of a $2$-chain $W$ then 
$m_\al(\ga)=\om(W)$ mod $\Z$.  Hence it contains the same information as 
the splitting $s^\Z$.
  For general $[\om]$, one can think of
a Hamiltonian structure as the marking defined by the splitting $s$;
cf. the discussion in the appendix of~\cite{LM}.}
Thus an  $M$-bundle $P\to B$  with structure group $\Ham^{s\Z}$ has such
monodromy functions on each fiber, but (just as in the case of 
the fiberwise
symplectic form) these do not need to be induced by a global monodromy function coming from a line bundle over $P$.  Such global structures
are called integral configurations in Gal--K\c edra~\cite{GalK}, where the problem of classifying them is discussed.
\end{rmk}

\begin{rmk}\labell{rmk:uniq}\rm {\bf (Issues of uniqueness)}
(i)    Because we are interested in the algebraic and geometric 
properties of the symplectomorphism group we restricted ourselves 
above to the case when $\Hh$ 
is a subgroup of $\Symp$.  However, 
from a homotopy theoretic point of 
view it would be more natural to look for a group $\Kk$ that 
classifies pairs consisting of a symplectic $M$-bundle $\pi:P\to B$
together with an extension $\Ta\in H^2(P;\R)$ of the fiberwise 
symplectic class $[\om]$.  
Here we should either 
normalize $\Ta$ by requiring  $\pi_!(\Ta^{n+1}) = 0$ (where $\pi_!$ 
denotes integration over the fiber) or consider 
 $\Ta$ to be well defined modulo elements in $\pi^*H^2(B)$. 
 Then the homotopy class of $B\Kk$ would be well defined and 
there would be a forgetful map $\psi:
B\Kk\to B\Symp$ which is well defined 
up to homotopy (assuming that we are working in the category of 
spaces with the homotopy type of a CW complex).  In general, $\psi$
 would not be a homotopy equivalence since the extension class $\Ta\in H^2(P;\R)$ could vary by an element in $H^1(B; H^1(M;\R))$. Further,  in this scenario, $\Kk$ need 
not be a subgroup of $\Symp$. (Cf. the discussion in
Lalonde--McDuff~\cite{LM} of the classification of Hamiltonian
structures.)
\sss
    
\NI
(ii) If we insist that $\Kk$ be a  subgroup of $\Symp$ then there are 
several
possible notions of equivalence, the most natural
of which is perhaps given by
conjugation by an element in $\Symp_0$. 
With this definition equivalent groups would be isomorphic.
We show in \S\ref{sec:Fs}  that the  groups $\Ham^s$ are equivalent 
in this sense, though when $H_1(M;\Z)$ has torsion the integer 
versions $\Ham^{s\Z}$ may not be.
It is also not clear whether any two 
groups $\Hh_1, \Hh_2$ that intersect each component of $\Symp$ and 
satisfy $\Hh_1\cap \Symp_0 = \Hh_2\cap \Symp_0=\Ham$
must be isomorphic 
as abstract groups, although any such group must be isomorphic to an 
extension of $\pi_0(\Symp)$ by $\Ham$.  Moreover, there is no 
immediate reason why they should
 be conjugate.  For example, suppose that 
the group $\pi_0(\Symp)$ is 
isomorphic to $\Z$, generated by the component $\Symp_{\al}$ of 
$\Symp$.
Then because $\Ham$ is a normal subgroup of $\Symp$ the subgroup 
$\Hh_{g}$ of $\Symp$, generated by $\Ham$ together with 
any element 
$g\in \Symp_{\al}$, intersects $\Symp_0$ in $\Ham$ and therefore
has the required properties.
Any two such groups $\Hh_{g_i}, i = 1,2,$
are isomorphic, though they are conjugate only if there is $h\in 
\Symp$ such that $g_1 hg_2^{-1}h^{-1}\in \Ham$.
 On the other hand, because $g_1$ and $g_2$
can be joined by an isotopy, there is a smooth family of
injective group homomorphisms $\iota_t:\Hh_{g_1}\to \Symp$, $t\in 
[1,2]$, that starts with the inclusion and ends with an isomorphism 
onto $\Hh_{g_2}$. Thus the homotopy properties of the inclusions
$\Hh_{g_i}\to \Symp$ are the same.  
\sss

\NI
(iii) Instead of looking for
 subgroups of $\Symp$ with the extension property one could look 
for {\it covering} groups  $\Hh\to \Symp$ with this property. 
Notice that  if $\La$ is a 
discrete subgroup of an abelian topological group $\Aa$ and if the 
continuous  crossed 
homomorphism $F:G\to \Aa/\La$ extends the composite $F_0:G_0\to 
\Aa\to 
\Aa/\La$, where $f:G_0\to \Aa$ is a homomorphism defined on the 
identity 
component of $G_0$, then the fiber product 
$$
\Tilde G: = \bigl\{(g,a)\in G\times \Aa\,|\, F(g) = a+\La\bigr\}
$$
of $G$ and $\Aa$ over $\Aa/\La$ 
is a covering group of $G$ that contains a copy of $G_0$, namely the 
graph of $f$.  Moreover, 
the obvious projection $\Tilde G\to \Aa$ lifts $F_0$.    
This approach is particularly relevant 
    in the integral case mentioned in Remark~\ref{rmk:hsz} above, as 
well
    as the cohomologically symplectic case, where the analog of 
    the Hamiltonian group is already a covering group of $\Diff_0$.
    For further discussion see \S\ref{sec:cHam} and 
    Gal--K\c edra~\cite{GalK}.
\end{rmk}

%%%%%%%%%%%%%%%%%%%%%%%%%%%%%%%%%%%%%%%%%%%%%%%%%%%%%%%%%%%%%%%%%%%%%%%%

\section{Definition and Properties of $\wF_s$}\labell{sec:Fs}
%%%%%%%%%%%%%%%%%%%%%%%%%%%%%%%%%%%%%%%%%%%%%%%%%%%%%%%%%%%%%%%%%%%%%%%%

Define $\Poz$ (resp. $\Po: = \Poq$) to be the set of values 
taken by $[\om]$ on the elements of $H_2(M;\Z)$
(resp. $H_2(M;\Q)$).  
To define $\wF_s$ we follow a suggestion of
 Polterovich (explained in Lalonde--McDuff~\cite{LM}).
Define the homology group 
 $$
 SH_1(M, \om;\Z)
 $$
to be the quotient of the  space of 
integral $1$-cycles in $M$ by the image under the boundary map $\p$ 
of  the integral
$2$-chains with zero symplectic area.  Then there is a projection 
$\pi_{\Z}: SH_1(M, \om;\Z)\to  H_1(M;\Z)$ and we  set
$$
  SH_1(M, \om):= SH_1(M, \om;\Z)\otimes \Q.
$$
We shall consider $SH_1(M, \om)$ and $\Poq$ as
$\Q$-vector spaces. 
Given a loop (or integral 
 $1$-cycle) $\ell$ in $M$ we denote its image  in 
 $H_1(M;\Z)$ or $H_1(M;\Q)$ by $[\ell]$  and its 
 image in $SH_1(M,\om;\Z)$ or $SH_1(M,\om)$ by $\langle\ell\rangle.$
We usually work over the rationals and shall omit the label $\Q$ 
unless there 
is a possibility of confusion.

\begin{lemma}\labell{le:split0}  There are split exact
sequences 
\begin{equation}\labell{eq:esz}
0 \to \R/\Poz\to SH_1(M,\om;\Z) \stackrel{\pi_{\Z}} \to H_1(M;\Z) \to 
0,
\end{equation}
and
\begin{equation}\labell{eq:es}
0 \to \R/\Poq\to SH_1(M,\om) \stackrel{\pi}\to H_1(M;\Q) \to 0.
\end{equation}
\end{lemma}
\begin{proof} Choose a continuous family of   
integral $2$-chains  $f_t: D\to M$ for $t\in \R$ with 
 $\int_{D}f_t^*\om = t$.  If  $\ga_t: = f_t|_{\p D}$ denotes the 
 boundary of $f_t$, then the 
 elements 
$$
\langle \ga_t\rangle , \quad t\in \R,
$$
generate the 
kernel of the projection $\pi_{\Z}: SH_1(M,\om;\Z) \to H_1(M,\Z)$.  
Moreover they represent 
different classes in $SH_1(M,\om;\Z)$ if and only if $t-t'\notin 
\Poz$.
Hence the sequence
$$
0 \to \R/\Poz\to SH_1(M,\om;\Z)  \stackrel{\pi_{\Z}}\to H_1(M;\Z) \to 
0,
$$
is exact.  
To see that it splits, we just need to check that each element 
$\la = [\ell]$ of 
finite order  $N$ in $H_1(M;\Z)$  is the image of some element of 
order $N$ in $SH_1(M;\Z)$.  But if  $W$ is an
integral $2$-chain such that $\p W = N\ell$ and if $\mu: = \int_W\om$ 
then
$$
N\bigl(\langle \ell\rangle - \langle \ga_{\mu/N}\rangle\bigr) = 
0\;\mbox{ 
and }\;
 \pi\bigl(\langle \ell\rangle - \langle \ga_{\mu/N}\rangle\bigr) = 
[\ell].
$$
In fact every element of order $N$ in the coset
$\pi_\Z^{-1}([\ell])$ has the form
$\langle \ell\rangle - \langle 
\ga_{\nu}\rangle$ where $N\nu\in \mu + \Poz$.
 The proof for~(\ref{eq:es}) is  similar.
\end{proof}

We explain in \S\ref{ss:Z} a natural way to understand
 splittings  of $\pi_\Z$ in the case when $\Poz = \Z$:
 cf. Definition~\ref{def:scanon}.  Note also 
 that in the previous lemma there is no need for $\om$ to be 
 nondegenerate; it suffices for it to be closed.
 However if it were an arbitrary closed form
it would not have many isometries, and so the next 
lemmas would have little interest.

\begin{lemma}\labell{le:split1} 
The group $\Symp(M,\om)$ acts  on $SH_1(M,\om;\Z)$ and
$SH_1(M,\om)$.  The induced action of $\Symp_0$
on the set of splittings of $\pi$ is 
transitive.  When $H_1(M;\Z)$ has no torsion  $\Symp$ also acts 
transitively
on the splittings of $\pi_{\Z}$.
\end{lemma}
\begin{proof} 
Again, we shall work with the sequence over $\Z$.
    The group $\Symp(M,\om)$ acts 
on these spaces  because it preserves $\om$.  
To prove  the transitivity statement, note first that any 
splitting $s$  of $\pi_{\Z}$ has the form $s\la_i = 
\langle\ell_i\rangle$
where $\ell_1,\dots, \ell_k$
are loops (i.e. integral $1$-cycles)
 in $M$  that project to the basis 
$\la_1,\dots,\la_k$ of 
$H_1(M;\Z)$.  Suppose given two such 
 splittings $s,s'$  corresponding to different sets $L, L'$ of 
 representing
 $1$-cycles for the $\la_i$.
Suppose also that  $\dim M > 2$.  Since Hamiltonian isotopies 
have zero flux, we may move the loops in $L$ and $L'$ by such 
isotopies,
without affecting their images in $SH_1(M,\om;\Z)$ and
so that no two intersect.
Now choose $T_1,\dots,T_k\in \R^+$ such that 
$$
\langle\ell_i'\rangle = \langle\ell_i\rangle +
\langle\ga_{T_i}\rangle,\quad 1\le i\le k,
$$
where the $\ga_t$ are as in Lemma~\ref{le:split0}.
For each $i$ there is a symplectic isotopy $h_{i,t}$ such 
that for all $t\in [0,T_i]$,
$$
h_{i,t}|_{\ell_j} = \id,\; j< i,\quad
h_{i,t}|_{h_{j,T_j}\ell_j} = \id, \; j> i,\quad
\int_{W_i}\om  = T_i,
$$
where $W_i: = \cup_{0\le t\le T_i} h_{i,t}(\ell_i)$.
(Take the $h_{i,t}$ to be generated by  closed $1$-forms $\al_i$ that 
vanish near the appropriate loops and are such that
$\int_{\ell_i}\al_i \ne 0$.  Here we are using the fact that 
$[\ell_i]$ is not a torsion class.)  
Then $h: = h_{1,T_1}\circ\dots\circ h_{k,T_k}$ takes $s$ to $s'$.

To extend this argument to the case $\dim M = 2$, it is convenient to 
describe the splitting by its effect on a standard basis  $\la_i$ of 
$H_1(M;\Z)$.  Thus we may assume that  $\ell_i$ and $\ell_j$ 
are disjoint unless $(i,j) = (2k-1,2k)$ in 
which case they intersect  in a single point.  If $s_0$ is the 
splitting defined by these loops, it suffices to show that for any 
numbers $T_i$ there are representatives $\ell_i'$ for the $[\ell_i]$ 
such that for each $i$ there is a cylinder of area $T_i$ with boundary
$\ell_i' -\ell_i$.  One achieves this by first isotoping
the $\ell_i$ for $i$ odd (fixing the other loops), and 
then adjusting the $\ell_i$ for even $i$.  
\end{proof}

 Choose a splitting $s$ for $\pi_\Z$. If $h\in \Symp$ and $\la\in 
 H_1(M;\Z)$, then the element $h_*(s\la) - s(h_*\la)$ lies in the 
 kernel of $\pi_\Z: SH_1(M,\om)\to H_1(M;\Z)$ and one can define
 a map
$$
 \wF_s^{\Z}: \Symp(M)\to \Aa^{\Z}: = \Hom\bigl(H_1(M;\Z),\R/\Poz\bigr)
$$
 by setting
 \begin{equation}\labell{eq:fsz}
 \wF_s^{\Z}(h)(\la):= 
h_*(s\la) - s(h_*\la)
 \in \R/\Poz,\quad \la \in H_1(M;\Z).
 \end{equation}
 Explicitly, if we
    denote by $\ov\la$ the image $s(\la)$ 
of $\la \in H_1(M)$, then 
\begin{equation}\labell{eq:Fs}
\wF_s(h)(\la) = a\bigl(h\ov\la-\ov{h\la}\bigr),
\end{equation}
where $a(\langle\ell'\rangle - \langle\ell\rangle)$ is the 
symplectic area of any cycle with boundary $\ell' - \ell$.
Similarly, for each splitting $s$ of (\ref{eq:es})
we define
 $$
 \wF_s: \Symp(M)\to \Aa: = \Hom\bigl(H_1(M),\R/\Poq\bigr) 
 = H^1(M;\R)/H^1(M;\Poq)
 $$
by 
 $$
 \wF_s(h)(\la):= 
h_*(s\la) - s(h_*\la)
 \in \R/\Poq,\quad \la \in H_1(M;\Q).
$$

\begin{prop}\labell{prop:wf} \NI
{\rm (i)}\,\, 
$
 \wF_s^{\Z}$ is a crossed homomorphism
that equals the composite
$$
\Symp_0\stackrel {\Flux}\to H^1(M;\R)/\Ga \to %H^1(M;\R/\Poz)\to
\Aa^{\Z}: = \Hom\bigl(H_1(M;\Z),\R/\Poz\bigr)
$$
on $\Symp_0$. Moreover the class 
$$
[\wF^{\Z}]: = [\wF_s^{\Z}] \in H_{cEM}^1\bigl(\Symp, \Aa^{\Z}\bigr)
$$
is independent of the choice of $s$.\sss

\NI
{\rm (ii)} The analogous statements hold for $\wF_s$. 
\end{prop}

\begin{proof} 
 $\wF_s^{\Z}$ is a crossed homomorphism because for all $g,h\in \Symp$
\begin{eqnarray*}
\wF_s^{\Z}(gh)(\la) & = & a\bigl(gh\ov\la-\ov{gh\la}\bigr)\\
& = &  a\bigl(gh\ov\la- g\ov{h\la}\bigr) + 
a\bigl(g\ov{h\la}-\ov{gh\la}\bigr)\\
& = & a\bigl(h\ov\la- \ov{h\la}\bigr) + 
a\bigl(g\ov{h\la}-\ov{gh\la}\bigr)\\
& = & \wF_s^{\Z}(h)(\la) + \wF_s^{\Z}(g)(h\la)\\
& = & \wF_s^{\Z}(h)(\la) + h^*\wF_s^{\Z}(g)(\la).
\end{eqnarray*}
The rest of the first statement in (i)
is immediate from the definition. 

To prove the second statement in (i) observe that two
choices of splitting $s,s'$ differ by the element
$\al\in \Aa^\Z: = \Hom\bigl(H_1(M;\Z),\R/\Poz\bigr)$ given by
$$
\al(\la):= s'(\la) - s(\la) \in \R/\Poz,\quad \la\in H_1(M;\Z). 
$$
It follows easily that
\begin{equation}\labell{eq:ss'}
 \wF_{s'}^{\Z}(h) - \wF_s^{\Z}(h) = \al - h^*\al,
\end{equation}
and  so is a coboundary.

The proof of  (ii) is similar.
\end{proof}

\begin{defn} Given a splitting $s: H_1(M;\Q)\to SH_1(M,\om)$
 we define the {\bf enlarged Hamiltonian group}
 $\Hams\Mo$ to be the kernel of $\wF_s$.  Similarly,
we define $\Ham^{s\Z}\Mo$ to be the kernel of the integral crossed 
homomorphism $\wF_s^{\Z}$.
 \end{defn}

\begin{lemma}\labell{le:conj}  Let $s, s' $ be two splittings and
 define 
 $\Symp_H$  to be the subgroup  of $\Symp$ that acts trivially on 
rational  homology.  \sss

 \NI
 {\rm (i)}  $\Symp_H\cap\, \Hams = \Symp_H\cap\, \Ham^{s'}$.\smallskip

  \NI
 {\rm (ii)} The map $\pi_{0}(\Ham^s) \to \pi_0(\Symp) $ is 
 surjective.
 
  \NI
 {\rm (iii)} The subgroups $\Hams$ and $\Ham^{s'}$ are conjugate 
in $\Symp$ by an element in $\Symp_0$.\sss

\NI
{\rm (iv)} When topologized  as a subspace of $\Symp$, 
the path component of $\Hams$  containing the identity element is 
$\Ham\Mo$.
\end{lemma}
\begin{proof}{}  (i) is an immediate consequence of the identity
(\ref{eq:ss'}).  
  (iii) follows
    from the fact that $\Symp_0$ 
acts transitively on the set of splittings and the
description of $\Hams$ as the subgroup of $\Symp$ whose action on
$SH_1\Mo$ preserves the image of $s$.  To prove (ii), we must
show that any element $h\in \Symp$ is homotopic to an element in 
$\Ham^s$.  This holds because the splittings $s' = h_*(s)$ and $s$ 
are  conjugate by an element in $\Symp_0$.
To prove  (iv) consider a 
continuous path 
$h_t\in \Symp$ that starts at the identity and is such that 
$\wF_s(h_t)(\la) = 0 \in \R/\Po$ for all $t$.  By 
Proposition~\ref{prop:wf}, 
the path $t\mapsto \wF_s(h_t)(\la)\in \R/\Po$ has the continuous lift
$t\mapsto \Flux(h_t)(\la)\in \R$.  Since $\Po$ is 
totally disconnected  this lift must be identically zero; in other 
words the path $h_t$ is a Hamiltonian isotopy.
 \end{proof}

Part (iv) of Lemma~\ref{le:conj} holds for  the group
$\Ham^{s\Z}$, and (i) holds if one replaces $\Symp_H$ by the group 
that acts trivially on $H_1(M;\Z)$.  However,
one must take care with the other two statements.
For further details see \S\ref{ss:Z}.

\begin{rmk}\labell{rmk:top}\rm
    {\bf (Topologies on $\Hams$ and $\R/\Po$.)}\,\,
The intersection $\Hams\cap \Symp_0$ is disconnected.  In 
fact it is everywhere dense in $\Symp_0$. Hence the subspace 
topology $\tau_s$
on $\Hams$ is rather counterintuitive and it is better to give  
$\Hams$ a finer topology in which its path components are
closed.  Therefore, although we give the group $\Symp$ the usual 
$C^{\infty}$-topology (which is the subspace topology it inherits 
from 
the diffeomorphism group), we give $\Hams$ the topology 
$\tau_c$ that it inherits from the {\it Hamiltonian} topology on
$\Symp$.
Then the identity map
$(\Hams,\tau_c)\to (\Hams,\tau_s)$ is continuous and is a weak 
homotopy equivalence.   Thus this change in topology does not affect 
the homotopy or (co)homology of the space.  

Correspondingly we shall 
always think 
of $\Po$ as a discrete group.
Further we think of 
quotients such as $\R/\Po$ (or $\Aa$) as quasitopological spaces\footnote
{
This meaning of the word \lq\lq quasitopological" seems to be obsolete,
though it is very convenient.}
 i.e.
 we specify which maps $f:X\to \R/\Po$ are continuous, where $X$ is a 
 finite simplicial complex.  This
gives enough structure for $\R/\Po$ to have well defined  homotopy 
groups and hence a well defined weak homotopy type.  In the present situation
we say that $f$ is continuous iff  $X$ has a triangulation $X'$ such 
that the restriction of $f$ to each 
simplex in $X'$ has a continuous 
lift to $\R$.  Hence 
$$
\pi_1(\R/\Po) \cong \Po,\qquad \pi_j(\R/\Po) = 0, \; j>1.
$$
Thus $\R/\Po$ is 
(weakly homotopic to) the Eilenberg--MacLane space $K(\Po,1)$.
In some situations (such as Lemma~\ref{le:obst} below) it is useful to replace the quasitopological space  $\R/\Po$  by its homotopy quotient
$\R/\!\!/\Po: = \R\times _\Po E\Po$. Here $E\Po$ denotes a contractible space 
on which the group $\Po$ acts freely and the notation $/\!\!/$ is 
taken from Segal~\cite{Seg}.  Since there is a homotopy fibration 
$\Po\to E\Po\to \R/\!\!/\Po$, $\R/\!\!/\Po$ is a $K(\Po,1)$.  Hence the map 
$\R/\!\!/\Po\to \R/\Po$ induced by collapsing $E\Po$ to a point is a weak
 homotopy equivalence.
  %Another way to deal with this technical problem --- that also 
%arises when one deals with spaces of germs --- is to replace 
%$\R/\Po$ by an appropriate semisimplicial complex: cf. May~\cite{May}.  But then to make sense of the homomorphism one has to 
%replace all spaces and groups by their semisimplicial analogs. 
%The only 
%place where we use these concepts is in the proof of
%Lemma~\ref{le:obst} and we explain there how to deal with these
%technicalities.
\QED\end{rmk}

 \begin{rmk}\labell{rmk:comm}\rm
If $g\in \Symp$ and
 $h\in \Symp_0$ then it is easy to check that
 $\Flux (g^{-1}hg) = g^*(\Flux h)$.
 Hence,  if $\Symp_H$ denotes the subgroup of $\Symp$ acting trivially
 on $H^1(M;\R)$ then
     $$
     [\Symp_0, \Symp_H] = \Ham. 
     $$
     On the other hand 
     $
     [\Symp_0, \Symp] = \Ham
     $
     only if 
     $\Symp = \Symp_H$. Hence when $\Symp \ne \Symp_H$ the flux 
     homomorphism does not extend to a homomorphism $\Symp\to 
H^1(M)/\Ga$.

     There is another relevant subgroup, namely $\Symp_{H^{\Z}}$, 
     consisting of elements that act trivially on $H_1(M;\Z)$.  Note
     that 
     $[\Symp_{H^{\Z}}, \Symp_{H^{\Z}}] \cap \Symp_0$ 
      lies in $\Ham^{s\Z}$ because $\wF_s^{\Z}$ restricts to a 
     homomorphism on $\Symp_{H^{\Z}}$ and so vanishes on the 
commutator 
     subgroup
     $[\Symp_{H^{\Z}}, \Symp_{H^{\Z}}]$.  But this is the best we can 
     say; in particular, it is not clear whether
     $[\Symp_H, \Symp_H] \cap \Symp_0$ 
     must always equal $\Ham$.  
\end{rmk}

\begin{lemma} \labell{le:commut} The following statements are 
equivalent.\smallskip
     
\NI
{\rm (i)}  $[\Symp_H, \Symp_H] \cap \Symp_0 = \Ham$;\smallskip
     
\NI{\rm (ii)}  For every product of commutators 
$y: = [u_1,u_2]\dots[u_{2p-1},u_{2p}]$, $u_i\in \Symp_H$, that lies in $\Symp_0$, there are elements $g_1,\dots,g_{2p}\in \Symp_0$ such that 
$$
f: = [u_1g_1,u_2g_2]\dots[u_{2p-1}g_{2p-1},u_{2p}g_{2p}] \in \Ham.
$$

     \NI
 {\rm (iii)} The flux homomorphism $\Flux:\Symp_0\to H^1(M;\R)/\Ga$ 
extends to
 a continuous homomorphism $F:\Symp_H\to H^1(M;\R)/\Ga$.
\end{lemma}
 \begin{proof}  Clearly (iii) implies (i), which in turn
  implies (ii).  To see that (ii) implies (i),
 note the identity
 $$
 [ug,vh] = g^u h^{uv} (g^{-1})^{uv} (h^{-1})^{uvu^{-1}}\,[u,v],
 $$
 where $g^a: = aga^{-1}$. It follows that
 $f y^{-1}$ may be written as a product
 of terms of the form 
 $g'_{2i-1}g'_{2i}(g''_{2i-1})^{-1}(g''_{2i})^{-1}$
 where $g_j'$ and $g_j''$ are conjugate to $g_j$ by products of the 
 $u_i$.   Since the $u_i$ lie in $\Symp_H$, $\Flux(g_j') =  
 \Flux(g_j'')= \Flux g_j$.  Hence $\Flux y = \Flux f = 0$, and $y\in \Ham$.
 
It remains to show that (i) 
 implies (iii).  As in the proof of Proposition~\ref{prop:Flux} given in \S1,
  it suffices 
 to find a section 
     $$
     \si:\pi_0(\Symp_H)\to \Symp_H,\quad \al\mapsto \si_{\al}\in 
\Symp_{\al},
$$     
 such that
 \begin{equation}\labell{eq:commut1}
     \si_{\al\be}\si_{\be}^{-1}\si_{\al}^{-1}\in \Ham,\qquad 
\al,\be\in 
     \pi_0(\Symp).
\end{equation}
We first  define $\si$ on the commutator subgroup 
      $[\pi_0(\Symp_H),\pi_0(\Symp_H)]$.
When $\al $ lies in this group then the component $\Symp_\al$ 
contains 
     elements that are products of commutators.  We define 
$\si_{\al}$  to be
     such an element.  Then $\si_{\al}$ is well defined modulo an 
element 
     in $\Ham$ because
     $
     [\Symp_H,\Symp_H]\cap \Symp_0 = \Ham
     $
   by assumption.  Hence~(\ref{eq:commut1}) holds for these $\al$.
   Now we extend by hand, defining a lift on the 
   abelian group 
  $
  \pi_0(\Symp_H)/[\pi_0(\Symp_H),\pi_0(\Symp_H)].
  $
  This is easy to do on the free 
  part, and on the torsion part one uses the divisibility of 
$H^1(M;\R)/\Ga$.
  Note that $F$ is necessarily continuous since it is continuous on 
  $\Symp_0$.
    \end{proof}

%%%%%%%%%%%%%%%%%%%%%%%%%%%%%%%%%%%%%%%%%%%%%%%%%%%%%%%%%%%%%%%%%%%%%%%%
\section{Bundles with structural group $\Hams$}\labell{sec:Hams}
%%%%%%%%%%%%%%%%%%%%%%%%%%%%%%%%%%%%%%%%%%%%%%%%%%%%%%%%%%%%%%%%%%%%%%%%

This section contains the proofs of the main results about the group 
$\Ham^s$ and the obstruction class. In \S\ref{ss:ext} we give a 
simple proof
that $\Symp^{\ham}$ has the extension property 
(Proposition~\ref{prop:sympham}).  Because
$\Hams$ is geometrically defined, a similar argument shows that 
$\Hams$ 
has the extension property  when restricted to bundles $P\to B$ 
where $\pi_1(B)$ acts trivially on $H_1(M;\Q)$.  However, in the case when $[\om]$ is integral, the analogous statement for $\Ham^{s\Z}$ holds only under
additional hypotheses: see Proposition~\ref{prop:main}.

We start \S\ref{ss:obs} by defining the obstruction cocycle $\OM$ 
and proving Theorem~\ref{thm:obs}.  (As noted by Tsemo~\cite{Tse}, 
this is a special case of a more general
theory that can be nicely expressed in the language of gerbes.)  We
then prove Propositions~\ref{prop:HHam1},~\ref{prop:modext}
and~\ref{prop:ham}.  The section ends with a proof of the
stability result Proposition~\ref{prop:stab}.

%%%%%%%%%%%%%%%%%%%%%%%%%%%%%%%%%%%%%%%%%%%%%%%%%%%%%%%%%%%%%%%%%%%%%%%%
\subsection{Groups with the extension property}\labell{ss:ext}
%%%%%%%%%%%%%%%%%%%%%%%%%%%%%%%%%%%%%%%%%%%%%%%%%%%%%%%%%%%%%%%%%%%%%%%%

We  begin by proving Proposition~\ref{prop:sympham}
which states that the group $\Symp^{\ham}$ has the extension 
property.\MS

\NI
{\bf Proof of Proposition~\ref{prop:sympham}.}
Suppose first that a smooth $M$-bundle $P\to B$ has a closed 
connection 
form $\Om.$  Because the holonomy of the corresponding connection is 
Hamiltonian round all contractible loops, it defines a {\it 
continuous} 
map from the space of based loops in $B$ to the group 
$\Symp^{\ham}$.  
This deloops to a lift $B\to B\Symp^{\ham}$ of the classifying map 
for $P\to B$.  Therefore the classifying map of
any bundle with a closed connection form does lift to $B\Symp^{\ham}$.

Conversely, 
consider the universal $M$-bundle 
$$
M_{\Symp^{\ham}}\to B\Symp^{\ham}.
$$
It suffices 
to show that the fiberwise symplectic class $a=[\om]$ extends to a 
class $\Tilde a\in H^2(M_{\Symp^{\ham}};\R)$.  If not,
there is a map of a finite CW complex 
$X\to B\Symp^{\ham}$ such that the fiberwise symplectic class in the 
pullback bundle $M\to P_X\to X$ does not extend to $P_X$.
By embedding  $X$ in Euclidean space  and replacing it by a small 
 open  neighborhood, we 
 may assume that $X$ is a smooth 
 (open) manifold.
 Hence we may suppose that  
 $M\to P_X\to X$ is smooth.   Since the structural group is 
$\Symp^{\ham}$ 
 this bundle has a symplectic  connection 
with holonomy in $\Symp^{\ham}$.  The holonomy round contractible 
loops lies 
in the identity component of $\Symp^{\ham}$ and hence is 
Hamiltonian.  
Therefore the Guillemin--Lerman--Sternberg construction provides a 
{\it closed} connection form $\tau$ on $P_X$ that defines this 
connection: 
see~\cite[Thm~6.21]{MS}.  Since $[\tau]\in H^2(P_X)$  extends 
$[\om]$, this 
contradicts our initial assumption.
\QED\MS

\begin{cor}\labell{cor:hams}  Let $\Hh$ be any subgroup of $\Symp$
    whose identity component is contained in $\Ham$.      
    Consider the universal $M$-bundle
    $$
    M\to M_{\Hh}\to B\Hh.
    $$
Then the fiberwise symplectic class $a: = [\om]$  extends to $\Ta\in 
H^2(M_{\Hh};\R)$.
\end{cor}
\begin{proof}  The hypothesis on $\Hh$ implies that the inclusion 
$\Hh\to \Symp$ factors
continuously through $\Symp^{\ham}$.  Therefore the class 
$\Tilde a\in H^2(M_{\Symp^{\ham}};\R)$ constructed above pulls back 
to 
$H^2(M_{\Hh})$.
\end{proof}

\begin{lemma}\labell{le:integ} If $\Poz = \Z$ then the
universal $M$-bundle over $B\Ham^{s\Z}$ carries an extension
of $[\om]$ that takes integral values on all cycles lying over the $1$-skeleton of the base.  
%In other words, we may choose $\Ta$ to lie in the image
%of $H^2(M_{\Hh};\Z)$ in $H^2(M_{\Hh};\R)$.
\end{lemma}
\begin{proof}   The universal $M$-bundle over 
$B\Ham^{s\Z}$ carries a connection with holonomy in 
$\Ham^{s\Z}$.  Since this has Hamiltonian holonomy round closed 
loops, the GLS construction shows that it is given by
a closed connection form $\Om$.  We claim that $[\Om]$ takes integral values on all cycles in $\pi^{-1}(B_1)$, where $B_1$ is the $1$-skeleton of
$B\Ham^{s\Z}$. 
%there is $b\in H^2(B;\R)$ such that $[\Om]+\pi^*(b)$ is integral.
%
% that is normalized by the condition
%$$
%\int_{M}[\Om]^{n+1} = 0\in H^2(B\Ham^{s\Z}; \R).
%$$
%(Here one can either work with simplicial 
%de Rham forms or pull back over some smooth base $B$.)
%We claim that $[\Om]$ is integral.

Since $[\om]$ is assumed integral, we need only
check that $[\Om]$ takes integral values on
%the cycles in the fiber.  Hence we need only check that it is integral on 
%the it suffices to consider the pullback bundle $P\to \Si$ over 
%any map 
%$\phi:\Si\to B\Ham^{s\Z}$, where $\Si$ denotes a Riemann surface.
%There are three kinds of integral classes in  $H_2(P;\R)$:
%cycles lying entirely in the fiber, cycles that project nontrivially 
%to
%$H_2(\Si)$ and 
cycles $C(\ga,\de)$ formed as follows.  Suppose that
$\ga$ is a closed path in the base with holonomy $m_{\ga}:M\to M$
that fixes the class $\de\in H_1(M;\Z)$.  Choose a loop $\ell_{\de}$
in $M$ such that $\langle\ell_{\de}\rangle = s(\de)$, and define
$C(\ga,\de)$ to be the union of the cylinder $C'$ formed by the
parallel translation of $\ell_{\de}$ around $\ga$ with a chain $C''$
in $M$ with boundary $\ell_\de - m_{\ga}(\ell_\de)$.  
Now observe that since $\Om = 0 $
on $C'$ and $m_{\ga}\in \Ham^{s\Z}$, 
equation~(\ref{eq:Fs}) implies that
$$
\int_{C(\ga,\de)} \Om = \int_{C''}\om = -\wF_s^{\Z}(m_{\ga}) \in \Poz
\subset \Z.
$$
This completes the proof.
%
%Thus $[\Om]$ takes integral values on cycles of the first and third 
%kinds.  
%The class $[\Om]+ \pi^*(b)$ is represented in the $E_2$-term of
% the Leray--Serre spectral sequence for the cohomology of
%$P\to \Si$ by a sum whose  $E_2^{02}$ entry is $[\om]$
%and whose $E_2^{11}$ entry is the  unique element of $H^1(\Si, 
%\{H^1(M;\R)\})$
%that equals $[\Om]$ on the cycles $C(\ga,\de)$. Both these entries are 
%integral.  By suitable choice of $b\in H^2(B;\R)$ we can also arrange that 
%the $E_2^{20}$ entry is integral.  Hence result.
%We now claim that
%because of the normalization condition $[\Om]$ must be integral.
%One way to see this is to look at 
%the Leray--Serre spectral sequence for the cohomology of
%$P\to \Si$. It
%suffices to check that there is an integral element  in the 
%$E_2$-term that 
%maps down to a representative for $[\Om]$ in $E_{\infty}$. 
%But we can take the $E_2^{02}$ entry to be $[\om]$
%and the  $E_2^{11}$ entry to be the  unique element of $H^1(\Si, 
%\{H^1(M;\R)\})$
%%that equals $[\Om]$ on the cycles $C(\ga,\de)$. Hence the $E_2^{11}$ 
%entry is integral. Finally, the
%normalization condition implies that the $E_2^{20}$ entry must be
%zero.
\end{proof}

\begin{prop} \labell{prop:main}  {\rm (i)} Let $P\to B$ be a symplectic bundle 
 over a finite simplicial complex $B$ such that $\pi_1(B)$ acts 
 trivially on $H_1(M;\R)$.  Then $P$ has a closed connection form iff 
 the classifying map for $P\to B$ lifts to $B\Hams$.
 \smallskip
 
\NI  {\rm (ii)} Let $\pi:P\to B$ be a symplectic bundle 
 over a finite simplicial complex $B$ such that $\pi_1(B)$ acts 
 trivially on $H_1(M;\Z)$. Suppose further that either $H_2(B;\Z)$ is free or 
 $P\to B$ admits a section over its $3$-skeleton.  Then $P$ has a closed and integral 
 connection form iff 
 the classifying map for $P\to B$ lifts to $B\Ham^{s\Z}$.
\end{prop}
\begin{proof}
 Corollary~\ref{cor:hams} shows that every $\Hams$-bundle has a 
closed 
 connection form.  Conversely, suppose that  $P\to B$ has a 
 closed connection form.  Then the restriction map $H^2(P;\R)\to 
 H^2(M;\R)$ contains $[\om]$ in its image.  Because $\Q$ is a field, 
the 
 restriction map $H^2(P;\Poq)\to 
 H^2(M;\Poq)$ also contains $[\om]$ in its image.  
 Choose a class $a\in H^2(P;\Poq)$
 that extends $[\om]$.  Thurston's construction (cf. 
\cite[Thm~6.3]{MS}) 
 provides a closed extension $\Om$ in class $a$.  We claim that the 
 holonomy of $\Om$ round loops  $\ga$ in the base $B$ lies in  
 $\Hams$.  Granted this, one can use the local trivializations
 given by $\Om$ to reduce  the structural group to $\Hams$.
 
 To prove the claim,  observe that because 
  $\pi_1(B)$ acts 
 trivially on $H_1(M;\R)$ one can use the connection defined by $\Om$
 to construct for {\it each} loop $\ga$ in $B$ a
$2$-cycle $C(\ga,\de_i)$ as in Lemma~\ref{le:integ}, where $[\de_i]$
runs through a basis of $H_1(M;\Q)$.  Then, the $\Om$-holonomy
$m_{\Om}(\ga):M\to M$ round the loop $\ga$ in $B$ satisfies the
identity:
    $$
    \wF_s(m_{\Om}(\ga))(\de_i) = -\int_{C(\ga,\de_i)}\Om\in \Poq.
    $$
    Hence $m_{\Om}(\ga)\in\Hams$.   This completes the proof of (i).
    
    Now consider (ii).  If $P\to B$ has an integral closed connection form, then the argument given above shows that its structural group reduces to
    $\Ham^{s\Z}$.  Conversely, 
 Lemma~\ref{le:integ} shows that any bundle pulled back from $B\Ham^{s\Z}$
 has a closed connection form $\Om$ that takes integral values 
 on cycles lying over $B_1$. 
 In other words, $[\Om]$ takes integral values on the elements in 
 $\ker\pi_*: H_2(P;\Z)\to H_2(B;\Z)$.  Hence the homomorphism
 $H_2(P;\Z)\to \R/\Z$  induced by $[\Om]$ may be written as a composite $f\circ \pi_*$, 
 where $f: H_2(B;\Z) \to \R/\Z$.  It suffices to check that
$f$ lifts to a homomorphism $\be: 
H_2(B;\Z) \to \R$.  Thus we need $f$ to
vanish on the torsion elements of $H_2(B;\Z)$.  
This is obvious if $H_2(B;\Z)$ is free, while, if there is a section over the $3$-skeleton, $f$ vanishes on the torsion classes in $B$ 
because they lift to torsion classes in $P$.
\end{proof}

We show in the next section that part (i) of this proposition extends to arbitrary bundles.  However the integral case is more subtle.  
Example~\ref{ex:moreinteg} shows
that when $[\om]$ is integral
there might be manifolds $(M,\om)$ for which there is no group that classifies symplectic $M$-bundles with integral closed connection form.  
Moreover the following example shows that 
the obstruction to the existence of an integral extension of $[\om]$ does involve the $3$-skeleton of $B$. Take the universal $S^2$-bundle $P\to B\SO(3)$, and provide $S^2$ with a symplectic form $\om$ in the class that generates $H^2(S^2;\Z)$.   Then
the first Chern class $c: = c_1^{\rm Vert}$ of the vertical tangent bundle extends $2[\om]$.  But $\om$ has no integral extension;
$c/2$ is not integral (the restriction of $P$ over the $2$-skeleton is the one-point blow up of $\C P^2$)  and $[\om]$ has a unique extension because $H^2(B\SO(3);\Z) = 0$. 

%%%%%%%%%%%%%%%%%%%%%%%%%%%%%%%%%%%%%%%%%%%%%%%%%%%%%%%%%%%%%%%%%%%%%%%%
\subsection{The obstruction class}\labell{ss:obs}
%%%%%%%%%%%%%%%%%%%%%%%%%%%%%%%%%%%%%%%%%%%%%%%%%%%%%%%%%%%%%%%%%%%%%%%%

Let $\La$ be a
countable subgroup of $H^1(M;\R)$ that contains $\Ga$
and is  invariant under the action of $\pi_0(\Symp)$.
Denote by $\Aa_\La$ the quasitopological abelian group 
$
H^1(M;\R)/\La$. As explained in Remark~\ref{rmk:top},
$\Aa_\La$ is homotopy equivalent to a $K(\pi,1)$ with 
$\pi_1$ isomorphic to the free  abelian group $\La$.
 Let $\wF_\La:\Symp\to \Aa_\La$ be a 
crossed homomorphism  whose restriction to $\Symp_0$ factors through $\Flux.$  Denote its kernel by  $\Hh_\La$.
The next lemmas hold trivially when  $\La=0$, for in this case
 the inclusion $\Hh_\La\to \Symp$ is a homotopy equivalence.

\begin{lemma}\labell{le:obst} There is an obstruction class 
$$
\OML\in H^2(B\Symp; \La),
$$
such that the classifying map
$\phi:B\to B\Symp$ of a symplectic bundle  lifts to $B\Hh_\La$  iff 
$\phi^*(\OML) = 0$.  
% Moreover $\OM=0$ if
%there is a continuous crossed homomorphism $\Tilde 
%    F_s: \Symp\to H^1(M;\R) $ that lifts  $\wF_s$.
\end{lemma}
\begin{proof}\,  Suppose first that $\La$ is a discrete subgroup of
$H^1(M;\R)$ so that $\Aa_\La$ is a topological (rather than 
quasitopological) space.  Consider the fibration sequence
$$
\Hh_\La \to \Symp\stackrel{\wF_\La} \to \Aa_\La,
$$
that identifies  $\Aa_\La$ as  the homogeneous space  $\Symp/\Hh_\La$.  There is an associated homotopy fibration 
$$
\Aa_\La\to B\Hh_\La\to B\Symp.
$$
Because $\Aa_\La$ is a $K(\pi,1)$, there is a single obstruction to the 
existence of a 
section of this fibration, namely a class 
$\OML\in H^2(B\Symp;\pi_1(\Aa_\La)) = H^2(B\Symp; \La)$.

In the general case, we should replace $\Aa$ by the homotopy quotient
$\Symp/\!\!/\Hh_\La$ and consider the homotopy fibration
$$
\Symp/\!\!/\Hh_\La\to \Symp \backslash\!\!\backslash \Symp/\!\!/\Hh_\La\to
\Symp\backslash\!\!\backslash *,
$$
where $*$ denotes the one point space, and the notation $G
\backslash\!\!\backslash X/\!\!/ H$ denotes the realization of the 
topological category formed from the action of the group $G\times H$ on $X$, where $G$ acts by multiplication on the left 
and $H$ by multiplication on the right.
As explained in~\cite{Seg}, $\Symp \backslash\!\!\backslash\Symp$ is contractible, and so
 $
 \Symp \backslash\!\!\backslash \Symp/\!\!/\Hh_\La = 
(\Symp \backslash\!\!\backslash\Symp)/\!\!/\Hh_\La
$
 is homotopy 
equivalent to $B\Hh_\La$.  Hence the above fibration  is a model for 
$\Aa_\La\to B\Hh_\La\to B\Symp$, and the argument proceeds as before.
\end{proof}

We now suppose that $\La = H^1(M;\La')$, where $\La'$ is some countable subgroup of $\R$ containing the integral periods $\Poz$ of $[\om]$.  Thus $[\om]\in H^2(M;\La')$.
In the next lemma 
$$
d_2^{\La}: H^2(M;\La')\to H^2\bigl(B\Symp;H^1(M;\La')\bigr)
$$
 denotes the second differential in the 
Leray--Serre spectral sequence for the cohomology of $M_{\Symp}\to B\Symp$
with coefficients $\La'$.

\begin{lemma}\labell{le:OMobst}  $\OML=d_2^{\La}([\om])$.
%The image of  $\OML$ in 
% $H^2\bigl(B\Symp;H^1(M;\La')\bigr)$ is $d_2^{\La}([\om])$.
\end{lemma}

\begin{proof}  Give $B: = B\Symp$ a CW decomposition with one vertex $*$ and fix an identification of the fiber $M_*$ over this vertex with $M$.
We shall show that both $\OML$ and $d_2^{\La}([\om])$ may be interpreted as the first obstruction to defining a closed connection form $\Om$ over the $2$-skeleton $B_2$ of $B$
whose monodromy round the loops in $B_1$ lies in 
$\Hh_\La$.  Note that because every loop
in $B$ is homotopic to one in $B_1$ and because $\Ham\subset\Hh_\La$, any such connection does have monodromy in $\Hh_\La$ and hence does define
an $\Hh_\La$-structure.

Choose a closed extension $\Om'$ of $\om$ over the $1$-skeleton  $B_1$ of
 $B$.  Since $\Hh_\La$ intersects every component of $\Symp$ we may 
 suppose that the holonomy of $\Om'$ is contained in $\Hh_\La$.  Let $\al:
 D\to B_2$ be a $2$-cell attached via $\al:\p D\to B_1$.  Choose an 
 identification $\Psi$ of the pullback symplectic bundle $\pi_D: 
 \al^*(M_{\Symp})\to D$ with the product $D\times M\to D$ that extends the 
 given identification $M_*=M$.   Then $\Psi$ is well defined modulo 
 diffeomorphisms of $D
\times M$ of the form $(z,x)\mapsto (z,\psi_z(x))$ where $\psi_z\in \Symp_0, \psi_*=\id$.
Hence the induced identification of $H^2(\pi_D^{-1}(\p D))$ with 
$H^2(S^1\times M;\La) = H^2(M;\La) \oplus  H^1(S^1;\Z)\otimes H^1(M;\La)$
is independent of choices.
 
Consider the pullback $\Om'_D$ of $\Om'$ to $\p D\times M$. Let  $h_t, t\in [0,1]$  be the family of symplectomorphisms defining its characteristic flow, i.e. for each $x$ the paths $(t,h_t(x))$ lie in the null space of $\Om'_D$.  Then $\Om'_D$ 
represents the class $[\om] + [dt]\times a_D$ in 
$H^2(S^1\times M;\La)$ where $a_D = \Flux (\{h_t\})$.  Since $h_1\in \Hh_\La$ by construction, $a_D\in \La: = H^1(M;\La')$. Note that $\Om_D'$ extends to a closed form over $\pi_D^{-1}(D)$ iff $a_D = 0$.  Hence the cocycle $D\mapsto a_D$ represents the obstruction 
class $\OML$.  It is also clear from the interpretation of $d_2$ 
via zigzags given in Bott--Tu~\cite[Thm~14.14]{BT}, that this cocycle also represents
 $d_2^\La([\om])$.  The fact that we only consider forms $\Om'$ over $B_1$ 
 with monodromy in $\Hh_\La$ corresponds to the fact that we restrict the coefficients to $\La'$.
\end{proof}

\begin{cor}\labell{cor:main2}  Suppose that $\La'$ is divisible, i.e. is a module over $\Q$.
Then the group $\Hh_\La$ has the extension property.
\end{cor}
\begin{proof}  Since $\La'$ is divisible, $d_2^\La([\om]) = 0$ iff $d_2([\om]) = 0$ in the spectral sequence with coefficients $\R$.  Since 
$d_3^\La([\om])$ also vanishes by Lemma~\ref{le:2skel}, it follows that a symplectic bundle has a $\Hh_\La$-structure iff it has a closed connection form.
\end{proof}

In particular, this proves Corollary~\ref{main}.  We now prove 
the other statements in \S1.2.  When $\La' = \Po$, the rational period group of $[\om]$, we denote $\OML=: \OM$.  \MS

    \NI
 {\bf Proof of Theorem~\ref{thm:obs}.}\,\,  Part (i) is an 
immediate 
 consequence of  Lemma~\ref{le:OMobst}. Now
 suppose there is an extension 
 $\TF:\Symp\to H^1(M;\R)/\Ga$ of $\Flux$.  The proof of 
 Lemma~\ref{le:OMobst} 
  shows that the image of $\OMG\in H^2(B\Symp;\Ga)$ in 
 $H^2(B\Symp;H^1(M;\Po))$ is $d_2^{\om}([\om]) = \OM.$  Thus  $\OM$ takes values in $\Ga$.
 Hence it remains to prove the converse, i.e. that $\Flux$ extends if 
 $\OM$ takes values in $\Ga$.  
 
 As in the proof of Proposition~\ref{prop:Flux} given in \S1,
   it suffices
 to find a section 
     $$
     \si:\pi_0(\Symp)\to \Symp,\quad \al\mapsto \si_{\al}\in 
\Symp_{\al},
$$     
 such that 
 \begin{equation}\labell{eq:commut2}
    (\si_{\al\be})^{-1}\si_{\al}\si_{\be}\in \Ham,\qquad 
\al,\be\in 
     \pi_0(\Symp).
\end{equation}
To do this,  
consider the fibration sequence $\Aa\to B\Hams\stackrel{\pi}\to B\Symp$
of Lemma~\ref{le:obst}.
  By assumption the obstruction to the existence of a section $s: 
B\Symp\to B\Hams$ is an element of $H^2(B\Symp; \Ga)$ where $\Ga$ is identified with its image in $H_\Q=\pi_1(\Aa)$.  This means that for any 
compatible CW structures put on $B\Hams$ and $B\Symp$ one can choose a map 
$s:(B\Symp)_1\to (B\Hams)_1$ (where $B_1$ denotes the $1$-skeleton of $B$)
 so that $\pi\circ s \sim {\rm id}$ and so 
that the corresponding obstruction cocycle takes values in $\Ga$.  Choose a CW structure on $B\Symp$ with one vertex, and one $1$-cell $I\times g_\al$ for each component $\al\in \pi_0(\Symp)$. (This is possible because $\pi_0(\Symp)$ is 
  countable.)  Then for each pair $\al,\be$ in $\pi_0(\Symp)$ there is a $2$-cell $c_{\al,\be}$ with boundary 
  $(I\times g_{\al\be})^{-1}(I\times g_\al)(I\times g_\be)$. (There are other $2$-cells in $B\Symp$ coming from  the $1$-skeleton of a CW decomposition for $\Symp$, but these are irrelevant for the current argument.)  We define a CW structure on $B\Hams$ in a similar way.  
  Then the map $s$ takes each  $1$-cell $I\times g_\al$ in $B\Symp$ to a loop in $(B\Hams)_1$.  This loop
   is given by a word $w_\al$ in the elements of $\Hams$ that represents
   an element $h(w_\al)$ in $\Hams\cap\Symp_\al$.   
   
    The obstruction to extending $s$ over the $2$-cell $c_{\al,\be}$ is the homotopy class in $B\Hams$ of the loop corresponding to the word  
  $(w_{\al\be})^{-1}w_\al w_\be$.  
  This can be identified with the homotopy class 
  $$
  [h(w_{\al\be})^{-1}h(w_\al)h(w_\be)] \in \pi_0(\Hams\cap \Symp_0) \cong H^1(M;\Po)/\Ga
  $$
 of the element $h(w_{\al\be})^{-1}h(w_\al)h(w_\be)$.  To say the obstruction $\OM(c_{\al,\be})$ takes values in $\Ga$ 
  means that 
  this class lies in $\pi_0(\Hams) = \Ham$.  Hence  it is always  possible to extend $s$ over these $2$-cells (though it may not extend over the other 
  $2$-cells in $B\Symp$). Further if we define
  the section $\si: \pi_0(\Symp)\to \Symp$ by $\si(\al) : = h(w_\al)$
  then the identity~(\ref{eq:commut2}) holds.      This completes the proof.  \QED

\begin{rmk}\labell{rmk:KKM}\rm  We sketch an alternative way to prove this result 
based on the ideas in K\c edra--Kotschick--Morita~\cite{KKM}.  Suppose that $\OM$ takes values in $\Ga$ but does not vanish on the 
$2$-skeleton $B_2$ of $B\Symp$.  Then one can form a new bundle $P_2'\to B_2$ with vanishing obstruction class by appropriately twisting the given
 bundle over each $2$-cell $c$ in  $B_2$
for which $\OM(c)\ne 0$; for each such $c$ change  the bundle by
taking the connect sum with an appropriate bundle $Q\to S^2$.  By construction there is a closed connection form $\Om'$ on $P_2'$. 
Since we did not change the bundle over $B_1$, the construction 
in~\cite[Theorem~6]{KKM} gives a class 
$[F]\in H^1\bigl((B\Symp^\de)_1; H^1(M;\R)\bigr)$. (As described in~\S4.2,
this is induced by the difference 
$[\tom]-\io^*[\Om']\in H^2(\io^*(P_1);\R)$, where $\io: B\Symp^\de\to B\Symp$.)  This class $[F]$ is represented by an Eilenberg--MacLane cochain
on the group $\Symp^\de$ with values in $H^1(M;\R)$.  However, because we changed the bundle over the $2$-cells, it satisfies the cocycle condition only when projected to the quotient $H^1(M;\R)/\Ga$.  Hence it gives rise to a crossed homomorphism $\Symp\to H^1(M;\R)/\Ga$ that extends Flux.\end{rmk}

We next turn to the proof of Proposition~\ref{prop:HHam1} which we restate for the convenience of the reader. Since (ii) and (iii) are obviously equivalent, we omit (ii) here.

\begin{lemma} \labell{le:HHam1} The  following  conditions
are equivalent.
\sss

\NI
{\rm (i)} There is  an extension
$
\Tilde F:\Symp(M,\om)\to H^1(M;\R)/\Ga
$
of Flux.
\sss

\NI
{\rm (ii)}
 For every product of commutators $[u_1,u_2]\dots[u_{2p-1},u_{2p}]$, $u_i\in \Symp$, that lies in $\Symp_0$, there are elements $g_1,\dots,g_{2p}\in \Symp_0$ such that 
$$
[u_1g_1,u_2g_2]\dots[u_{2p-1}g_{2p-1},u_{2p}g_{2p}] \in \Ham.
$$

\NI
{\rm (iii)} For every  symplectic $M$-bundle $P\to \Si$ there is a 
    bundle $Q\to S^2$  such that the fiberwise connect sum $P\#Q\to 
    \Si\#S^2 = \Si$ has a closed connection form. 
\end{lemma}    
    \begin{proof}
Because $\OM\in H^2(B\Symp; H^1(M;\Po))$ where $\Po$ is divisible, $\Om$ lies in the image of $H^2(B\Symp; \Ga)$ iff  its pullback over every map $\phi:\Si\to B\Symp$ takes values in $\Ga$.
Since $\Ga$ is the image of $\pi_1(\Symp)$ under the Flux
homomorphism, for every element $\ga\in \Ga$ there is a corresponding 
 $M$-bundle $Q\to S^2$ which is classified by a map $\phi:S^2\to B\Symp$ 
 such that $\phi^*(\OM)([S^2]) = \ga$.  Further if $\phi:\Si\to B\Symp$ classifies the bundle $P\to \Si$ and $\psi:S^2\to B\Symp$ classifies
 $Q\to S^2$, the fiber connect sum $P\#Q\to \Si$ is classified by
 $\phi\vee\psi: \Si\vee S^2\to B\Symp$.  Therefore the equivalence of (i) and (iii) follows immediately from Theorem~\ref{thm:obs}(ii).
  
To see that
 (i) implies (ii) observe that  if we write
 $\Hh: = 
\ker \TF$,  then the quotient 
    group $\Hh/\Ham$ is isomorphic to $\pi_0(\Symp)$
     because $\Hh\cap \Symp_0 = \Ham$. Hence any representation in $\pi_0(\Symp)$ can be lifted to the subgroup $\Hh/\Ham$ of $\Symp/\Ham$.
     This immediately implies (ii).
     
     It remains to show that (ii) implies (iii), which we do by direct construction.  
Consider any symplectic 
bundle $\pi: P\to \Si_g$.
Decompose it into the union of a trivial 
bundle $M\times D\to D$ over a $2$-disc with a symplectic bundle $P'
\to \Si_g\smallsetminus D$.  Choose a flat symplectic connection $\Om'$ on $P'$ 
whose holonomy round the 
generators of $\pi_1(\Si_{g})$ is given by 
elements $u_i\in \Symp$.  Since the symplectic trivialization of $\pi$ over $D$ is determined up to a Hamiltonian isotopy, there is an identification of $\p P'$ with $M\times S^1$ that is well defined
up to a Hamiltonian loop $g_t\in \Ham, t\in \R/\Z$.  
Thus the holonomy of $\Om'$ round $\p P'$ gives rise to a path $h_t$ from the identity to $f: = [u_1,u_2]\cdots[u_{2g-1},u_{2g}]$ that is well defined modulo a Hamiltonian loop.   
Hence $\Flux(\{h_t\})\in H^1(M;\R)$ depends only on the choice of flat connection $\Om'$.   By (ii)
 we may homotop the $u_i$ (or equivalently choose $\Om'$) so that $f\in \Ham$.  Thus the obstruction
$\phi^*(\OM)([\Si]) = \Flux(\{h_t\})$ lies in $\Ga$.\end{proof}

\NI
{\bf Proof of Proposition~\ref{prop:modext}.}\,\,
Let $\Hh$ be a subgroup of $\Symp$ with identity component
 $\Ham$ and consider the image  $\im(\pi_0(\Hh))$ of
$\pi_0(\Hh)$  in $\pi_0(\Symp)$.
%the identity component 
%of $\Hh$ is $\Ham$.  
If $\Hh$ has the modified extension property then we must show that the intersection with $\im(\pi_0(\Hh))$ of 
every finitely generated subgroup $G$ of $\pi_0(\Symp)$ has finite index in $G$.  But otherwise there would be  a map of a finite wedge $V$ of 
circles into $B\bigl(\pi_0(\Symp)\bigr)$  such that no finite cover 
$\Tilde V$ 
of $V$ lifts into the image of $B\Hh$ in $B\bigl(\pi_0(\Symp)\bigr)$.
Since any bundle over 
a $1$-complex has a closed extension form, this contradicts our 
assumption on $\Hh$.

Conversely, assume that the cokernel of 
$\im(\pi_0(\Hh))$ in $\pi_0(\Symp)$ has the stated finiteness 
properties and let $\Symp_{\Hh}$ be the subgroup of $\Symp$ 
consisting of elements isotopic to $\Hh$.  
If $P\to B$ is classified by a map into $B\Hh$ then it has a closed connection form by Corollary~\ref{cor:hams}.  Therefore we just need to see that if $\phi:B\to 
B\Symp$ classifies a bundle with a closed extension form its pullback over some finite cover $\TB\to B$ lifts to $B\Hh$.

Observe first that
 the composite map $\pi_1(B)\to \pi_0(\Symp)/\pi_0(\Symp_\Hh)$
has finite image by hypothesis.  (Recall that we always assume $\pi_1(B)$ is finitely generated.)  Therefore we may replace $B$ by a 
finite cover such that the pullback bundle $\Tilde P\to \Tilde B$
is classified by a map $\Tilde\phi:\Tilde B\to B\Symp_{\Hh}$.  
     Set $\La\subset H^1(M;\R)$ 
equal to the (discrete) group $\Flux(\Hh\cap \Symp_0)$, and then 
define a crossed homomorphism
$$
F: \Symp_{\Hh}\to H^1(M;\R)/\La
$$
as in the proof of Proposition~\ref{prop:Flux} given in 
\S\ref{ss:extflux}.  Because the bundle
$\Tilde P\to \Tilde B$ has a closed connection form, the   
class $\Tilde\phi^*(d_2([\om])$ vanishes in $H^2(\Tilde B; H^1(M;\R))$.
The proof of Lemma~\ref{le:OMobst} shows that this class is the image of  $\OML\in H^2(\Tilde B; \La)$ under the map induced by the inclusion $\La\to H^1(M;\R)$.  By pulling back over a further cover if necessary, we may suppose that
$H^2(\Tilde B;\La)$ has no torsion. (Since $\pi_1(B)$ may  act nontrivially on  the coefficients $\La$, it is not enough to assume that
$H_1(\Tilde B;\Z)$ is free.) Hence this map is injective and
$\OML = 0$.  Thus $\Tilde\phi$ lifts to $B\Hh$ as required.\QED\MS

 \NI
{\bf Proof of Proposition~\ref{prop:ham}.}\,\,
The first claim is that $\Ham$ has the modified restricted extension 
property.
This is a corrected statement of the conclusions that one can draw
from the proof of Theorem 1.1 in~\cite{LM}.  The claim
also follows by arguing as in the proof of 
Proposition~\ref{prop:modext} using $\Flux$ instead of $\wF_s$; 
the argument can be greatly simplified 
because the group $\Symp_0$ acts trivially on the coefficients.  Here 
one should also note 
that if the cover $\Tilde B\to B$ is chosen so that $H_1(\Tilde B;\Z)$
has no torsion, then the boundary map $\delta: H^1(B; H_\R/H_{\Q})\to 
H^2(B;H_\Q)$ vanishes.

The second claim is that when $\Ga\ne 0$ the group
$\Ham$ does not have the extension property.
To see this choose a nonzero element $\be\in H^1(M;\R)\smallsetminus 
\Ga$ such that
$2\be \in \Ga$ and then choose $g\in \Symp_0$ with $\Flux(g) = \be$.
Consider the bundle $P\to \R P^2$ that is formed from the mapping 
torus bundle   
$$
    M_g: = M\times [0,1]/(x,1)\!\sim\!(gx,0)\;\longrightarrow \;\; S^1
$$
by attaching $M\times D^{2}$ by the map
$(x,e^{2\pi i t})\mapsto  (g_t(x), 2t)$ where $g_t$ is a path in 
$\Ham$
from the identity to $g_1:= g^{-2}\in \Ham$.  
The 
flat connection on $M_g$ pulls back to a connection with Hamiltonian monodromy round 
the boundary $M\times\p D$ and so extends  to
a closed connection form over the rest of $P$: cf. the proof of
Lemma~\ref{le:OMobst}.

We claim that this bundle has no
 Hamiltonian structure.  
To see this consider the classifying map $\phi:\R P^2\to B\Symp_0$.
Just as in the discussion before Lemma~\ref{le:obst} the 
homomorphism $\Flux:\Symp_0\to H_\R/\Ga$ defines an obstruction class
$\OMG\in H^2(B\Symp_0;\Ga)$ such that $\phi^*(\OMG)$ vanishes iff
the bundle $P\to\R P^2$ has a Hamiltonian structure.  Since 
$B(H_{\R}/\Ga)$ is a $K(\Ga,2)$, this class is 
the pullback to $B\Symp_0$ of the canonical generator of 
$H^2(K(\Ga,2);\Ga).$  
We claim that the composite map 
$$
\R P^2\to B\Symp_0 \to B(H_{\R}/\Ga) = K(\Ga,2)
$$ 
is not null homotopic.  Since $\R P^2$ is the $2$-skeleton of $\R 
P^{\infty} = K(\Z/2\Z; 1)$ and $K(\Ga,2)$ is homotopy equivalent to a 
product of copies of $ BS^{1}$, this assertion is equivalent to 
saying that
under the  map $B(\Z/2\Z)\to BS^1$ induced by the obvious inclusion
$\{\pm 1\}\to S^1$ the 
generator of $H^2(B S^1; \Z/2\Z)$ pulls back to a nonzero element of
$H^2(\R P^2;\Z/2\Z)$.  This is well known.  For a direct proof 
identify the $2$-skeleton of $BS^1 = \C P^\infty$ with the quotient 
$S^3/S^1$ and observe that the $\Z/2\Z$-equivariant map
$$
S^2\to S^3,\qquad (r,s,t)\mapsto \Bigl(r+is, \frac 1{\sqrt 
2}(t+it)\Bigr)\in 
S^3\subset\C^2,
$$
descends to a map $\R P^2\to S^2$ of (mod 2) degree $1$.
\QED\MS

%%%%%%%%%%%%%%%%%%%%%%%%%%%%%%%%%%%%%%%%%%%%%%%%%%%%%%%%%%%%%%%%%%%%%%%%
\subsection{Stability}\labell{ss:stab}
%%%%%%%%%%%%%%%%%%%%%%%%%%%%%%%%%%%%%%%%%%%%%%%%%%%%%%%%%%%%%%%%%%%%%%%%

We finally discuss the question of stability.\MS

\NI
{\bf Proof of Proposition~\ref{prop:stab}.}  
Let $\Nn(\om)$  be a path connected  
neighborhood of $\om$ in the space of forms annihilating $V_2(P)$ 
such that $P\to B$ has an $\om'$-symplectic structure for all 
$\om'\in \Nn(\om)$. Our aim is to shrink $\Nn(\om)$ so that
each such $\om'$ has a closed extension to $P$.  We claim that
for each map $\psi:\Si\to B$ of a Riemann surface into $B$ there is a 
homologous map $\psi':\Si'\to B$   such that the pullback bundle over 
$\Si'$ admits a closed extension of $\om'$, provided that $\om'$ is 
sufficiently close to $\om$.  Granted this, we may choose $\Nn(\om)$ 
so that this holds  for a finite set of $\psi_i$ that 
represent a set of generators for $H_2(M;\R)$ and all $\om'\in 
\Nn(\om)$.  It follows that the obstruction class $\Oo^M_{\om'}$ must 
vanish when pulled back to $B$, i.e.
that $(M,\om')\to P\to B$ has a closed connection form when $\om'\in 
\Nn(\om)$.

To prove the claim, consider a map $\psi:\Si: = \Si_p\to B$.
As in the proof of Lemma~\ref{le:HHam1}, we may assume that
the pullback bundle $\psi^*P\to \Si$ has a flat $\om$-symplectic 
connection
over $\Si\smallsetminus D^2$  whose holonomy $y$
 around the boundary of the disc $D^2$ 
may be expressed as:
$$
 y: = [u_1,u_2]\cdots[u_{2p-1},u_{2p}]\in \Ham(M,\om),\quad 
u_i\in\Symp(M,\om).
 $$
 Since $\Ham(M,\om)$ is a perfect group, we may, by increasing the 
genus of $\Si$ and choosing the flat connection on the extra handles 
to have Hamiltonian holonomy, assume that $y={\rm id}.$  
 By hypothesis on the deformation $\om'$, we can choose:
 \MS
 
 \NI$\bullet$ 
 a path $\om_t$ from $\om_0: = \om$ to $
  \om_1: = \om'$ in $\Nn([\om])$, and \smallskip
  
  \NI$\bullet$ $C^1$-small paths
  $g_{it}\in
\Diff_0(M)$ such that
 $u_ig_{it}\in \Symp(M,\om_t)$ for all $i$ and $t\in[0,1]$.
   (These may be constructed using the Moser method.) 
 \smallskip
 
 Since $y={\rm id}$  the smooth path
 $$
y_t = [u_1g_{1t},u_2g_{2t}]\cdots[u_{2p-1}g_{(2p-1)t},u_{2p}g_{(2p)t}]
 $$
is $C^1$-small and  lies in $\Symp_0(M,\om_t)$ for all $t$.  
 If we could arrange that $y_t\in  \Ham(M,\om_t)$ for each $t$ then 
 the connection could be extended to a Hamiltonian connection over 
the disc for all $t$ and the proof would be complete.
 
 We do this in two stages.  First we modify
  $\psi:\Si\to B$ to a map $\psi':\Si'\to B$ so that $y_t\in 
\Ham^s(M,\om_t)\cap\Symp_0(M,\om_t)$ for all $t$. 
To this end,  consider the subspace $V^1$  of $H^1(M;\R)$ 
generated by 
the elements  $u_i^*\al-\al$, where $i = 1,\dots,2p,$ and 
$\al$ runs through the elements of $ H^1(M;\R)$. If the elements 
$a,b\in \Symp(M,\om)$ are each homotopic to some $u_i, i=1,\dots,2p$ 
then
\begin{eqnarray*}
\wF^s([a,b]) & = & \wF^s(b^{-1}) + (b^{-1})^*\wF^s(a^{-1})  + 
(a^{-1}b^{-1})^*\wF^s(b) + (ba^{-1}b^{-1})^*\wF^s(a)\\ 
& = & -b^*\wF^s(b)  + (a^{-1}b^{-1})^*\wF^s(b) 
 - (b^{-1})^*a^*\wF^s(a)   + (ba^{-1}b^{-1})^*\wF^s(a),
\end{eqnarray*}
which is easily seen to lie in $V^1/\bigl(V^1\cap H^1(M;\Po)\bigr)$. 
Hence 
$$
\wF^s_{\om_t}(y_t) \in V^1/\bigl(V^1\cap H^1(M;{\mathcal 
P}_{\om_t})\bigr),\quad t\in [0,1].
$$
By compactness we can therefore find a finite collection of smooth 
families
$(v_{jt},\al_{jt}), j = 1,\dots,m$, such that
$$
\sum_{j=1}^m v_j^*\al_{jt}-\al_{jt} \in \wF^s_{\om_t}(y_t)
+ H^1(M;{\mathcal P}_{\om_t}),\quad t\in [0,1],
$$
where each $v_{jt}\in \Symp(M,\om_t)$ is a product of 
the elements $(u_ig_{it})^{\pm 1}, i=1,\dots,2p,$ and $\al_{jt}$ is a 
path in $H^1(M;\R)$ with initial point $\al_{j0}=0$.

For each pair $(v_{jt},\al_{jt})$ choose  a path
$\Tilde h_{jt}$  in $ \Tilde\Symp_0(M,\om_t)$ starting at $\id$ 
 such that $\Flux_{\om_t}\Tilde h_{jt}
= \al_{jt}$.  Then 
$$
\Flux_{\om_t}[v_{jt}^{-1},\Tilde h_{jt}^{-1}] = \al_{jt} - 
v_{jt}^*\al_{jt} = :\be_{jt}.
$$
Therefore there is a fibration $M\to Q_j\to \T^2$
that admits an   $\om_t$-symplectic structure for each $t$ 
and a flat connection over $\T^2\smallsetminus D^2$
whose boundary  holonomy has flux 
$\be_{jt}$.  As an $\om_t$-symplectic bundle, $Q_j$ is pulled back 
from a bundle over $S^1$ with holonomy $v_{jt}$. Our choice of 
$v_{jt}$ implies this bundle is a pullback of $P\to B$ by some map 
$\psi_j: S^1\to B$ that we can assume to be independent of $t$ (since 
the holonomy $\be_{jt}$ depends only on the homotopy class of 
$v_{jt}$.)
However the connection varies smoothly with $t$.
Therefore we can change the flux $\wF^s_{\om_t}(y_t)$ of the 
boundary $\om_t$-holonomy  of the chosen flat connection on
$f^*P\to \bigl(\Si\smallsetminus D^2\bigr)$ to 
$\wF^s_{\om_t}(y_t)+\be_{jt}$ by 
replacing $\psi:\Si\to B$ by the homologous map 
$$
\psi\# \psi_j: \Si\#\T^2\to B.
$$
 Repeating this process for $j = 1,\dots,m$ allows us to 
 perform the required modification.
 
 Therefore we have now arranged that 
 $y_t\in 
\Ham^s(M,\om_t)\cap\Symp_0(M,\om_t)$ for all $t$.
The  following continuity argument shows that in fact
 $y_t\in \Ham(M,\om_t)$ for all $t$, which 
 finishes the proof. 
 
 Observe that
  for each $[\ell]\in H_1(M)$ and $t\in [0,1]$, the number
$$
 \Phi(t)[\ell]: =   \int_{[0,t]\times S^1} 
\phi_{t,\ell}^*\om_t \;\in \;\R,\;\;\mbox {where }   \phi_{t,\ell}(r,s): = 
y_r(\ell(s)),
$$
projects to $\wF^s_{\om_t} (y_t)([\ell]) \in \R/
{\mathcal P}_{\om_t}$.  Since, by assumption, $y_t\in \Ham^s(M,\om_t)$,
we find that  $ \Phi(t)[\ell]\in {\mathcal P}_{\om_t}$ for all $t$.   
 But $\Phi(t)$ varies continuously with $t$ and $\Phi(0) = 0$.  Hence the
fact that $V_2(\om_t)\supseteq V_2(\om)$ implies that $\Phi(t) 
= 0$ for all $t$.  
  It remains to check that $\Phi(t)$ projects to $ 
\Flux_{\om_t}(y_t)\in H^1(M;\R)/\Ga_{\om_t}$ for all $t\in [0,1]$. 
But because 
 the $y_t$ are $C^1$-small, for each fixed $t\in [0,1]$ the path 
$\{y_r\}_{r\in [0,t]}$ 
 may be canonically homotoped to a path $\{y_{rt}'\}_{r\in [0,t]}$ in 
$\Ham(M,\om_t)$
by a Moser process that fixes its endpoints.  $\Flux_{\om_t}(y_t)$
 is given by integrating $\om_t$ over the corresponding chain
 $\phi_{t,\ell}':[0,t]\times S^1\to M$.  Since this is homotopic to
 $\phi_{t,\ell}:[0,t]\times S^1\to M$  mod boundary, we find that 
 $
\Phi(t)= \Flux_{\om_t}(y_t)$ mod $\Ga_{\om_t}$,
  as required.

\QED\MS

%%%%%%%%%%%%%%%%%%%%%%%%%%%%%%%%%%%%%%%%%%%%%%%%%%%%%%%%%%%%%%%%%%%%%%%%

\section{Further considerations}\labell{sec:furth}
%%%%%%%%%%%%%%%%%%%%%%%%%%%%%%%%%%%%%%%%%%%%%%%%%%%%%%%%%%%%%%%%%%%%%%%%

We begin by collecting together various observations about the groups 
$\Ham^{s\Z}$ in the case when $\Poz=\Z$.
We then  explain some situations in which 
$\wF_s$ lifts to a crossed homomorphism with values in 
$H^1(M;\R)/\Ga$.
This is followed by  a short discussion of $c$-Hamiltonian bundles
and covering groups. 

%%%%%%%%%%%%%%%%%%%%%%%%%%%%%%%%%%%%%%%%%%%%%%%%%%%%%%%%%%%%%%%%%%%%%%%%
\subsection{The integral case}\labell{ss:Z}
%%%%%%%%%%%%%%%%%%%%%%%%%%%%%%%%%%%%%%%%%%%%%%%%%%%%%%%%%%%%%%%%%%%%%%%%

We shall assume throughout this section that $\Poz=\Z$. 
Many (but not all) of our results have some analog in the general 
case.

We begin by considering
 the integral analog of Lemma~\ref{le:conj}.
If $H_1(M;\Z)$ is torsion free, then Lemma~\ref{le:split1} applies 
and 
the whole of this lemma extends.    But if this group has torsion 
then it is possible that (ii) does not hold.

 \begin{lemma}\labell{le:hsz0}  Suppose that $\Poz=\Z$ and set $\Tor: 
= \Tor(H_1(M;\Z))$. 
     Then:\sss
     
     \NI
     {\rm (i)}\, $\wF_s^{\Z}$ induces a 
     crossed homomorphism $C_s:\pi_0(\Symp) \to \Hom(\Tor,\R/\Z)$, 
whose 
     kernel equals the image of $\pi_0(\Ham^{s\Z})$ in 
     $\pi_0(\Symp)$.  \sss
     
\NI
      {\rm (ii)}\,\, The image $[C_s]$ of $C_s$ in $H^1_{\EM}(\Tor; 
     \R/\Z)$ is independent of the choice of $s$. In particular, 
     if
      $\Symp$ acts trivially on $\Tor$ then
      the kernel of $C_s$ is independent of the choice of splitting 
      $s$.\sss
      
      \NI
      {\rm (iii)}\,\, There is a splitting such that $C_s=0$ iff 
$[C_s]=0$.
     \end{lemma}
\begin{proof}
We saw in Lemma~\ref{le:conj} that if $[\ell]$ has 
order $N$ there are precisely $N$ distinct elements of order $N$
in the coset $\pi_{\Z}^{-1}([\ell])$, 
namely $\langle \ell\rangle - \langle 
\ga_{(i+\mu)/N}\rangle$ for $i=0,\dots,N-1$, where $\mu$ is 
the area of a chain $W$ that bounds 
$N\ell$.  Since $\Symp_0$ is a connected group it must act trivially 
on these elements.   Therefore, for each $g\in \Symp$
the restriction of $\wF_s^{\Z}(g)$ to the 
torsion elements in $H_1(M;\Z)$  depends only on the image of $g$ in 
$\pi_0(\Symp)$.  This shows that $C_s$ exists.  Its kernel obviously 
contains the image of
$\Ham^{s\Z}$.  To complete the proof of (i)  we must show that
if $\wF_s^{\Z}(g)$ vanishes on $\Tor$ then $g$ may be isotoped to an 
element in $\Ham^{s\Z}$. But this holds by the proof of 
Lemma~\ref{le:split1}. (Note that we may assume that $\dim M >2$ here 
since otherwise $\Tor = 0$.)

Statement (ii) holds by the argument of Proposition~\ref{prop:wf}: 
given 
a splitting $s$ of $\pi_\Z$ over $\Tor$,
any other splitting $s':\Tor\to \pi_{\Z}^{-1}(\Tor)$ has the form
$s+ \langle\ga_{\be([\ell])}\rangle$  where   $\be\in 
\Hom(\Tor;\R/\Z)$.   
(iii) is an immediate consequence of (ii).
\end{proof}

\begin{example}\labell{ex:moreinteg}\rm (i) We again assume that 
$\Poz=\Z$.
    Suppose that 
    for some $g\in \Symp$  
    there is a loop $\ell$ such that $[\ell]$ has order $N>1$ in
    $\Tor$ and there is a $2$-chain $W$ 
    with boundary $g(\ell)-\ell$ and area $i/N$, where $0<i<N$. 
    Then $C_s(g)\ne 0$ for all  splittings $s$, and there 
 is no splitting such that $g$ is 
 isotopic to an element in $\Ham^{s\Z}$.  
Note that
    the corresponding mapping torus bundle 
$$
    M_g: = M\times [0,1]/(x,1)\!\sim\!(gx,0)\;\longrightarrow \;\; S^1
$$
    has a closed connection form but not one that is 
integral.\footnote
{This example is very similar to that in the proof of 
Proposition~\ref{prop:ham}. But now $[\ell]
\in H_1(M;\Z)$ is a torsion element,
and $g$ is not isotopic to the identity though it acts trivially on 
$[\ell]$.}
Equivalently, 
 $g$ does not fix any integral lift $\tau\in H^2(M;\Z)$ of $[\om]$.\MS

\NI (ii)
    The following yet more intriguing situation cannot be ruled 
out  in any obvious way.  Suppose  that 
$\Tor = \Z/2\oplus \Z/2$ is generated by the elements $[\ell]$ and 
$[\ell']$ which are interchanged by two symplectomorphisms $h_1, h_2$.
Suppose further that 
$$
(h_1)_*\langle\ell\rangle\ne (h_2)_*\langle\ell\rangle,\qquad
(h_1^2)_*\langle\ell\rangle= 
(h_2^2)_*\langle\ell\rangle=\langle\ell\rangle.
$$
Then $g: = h_1h_2$ fixes $[\ell]$ but acts nontrivially on $
\langle\ell\rangle$ and so has the properties assumed in (i) above.
Now consider the
splittings $s_i$ defined  by
$$
s_i[\ell] = \langle\ell\rangle,\qquad s_i[\ell'] = \langle 
h_i(\ell)\rangle,\qquad i=1,2.
$$
Then $h_i\in \Ham^{s_i\Z}$ by construction.
The corresponding mapping tori $M_{h_{i}}$ have 
$\Ham^{s_i\Z}$-structures and  so
 Lemma~\ref{le:integ} implies that each supports a closed integral 
connection form.   But their fiber connect sum $P\to V: = S^1\vee S^1$
does not, since one of its pullbacks is the bundle $M_g\to S^1$ 
considered 
in (i).
Indeed, the torsion in $H_1(M;\Z)$ creates new terms in  $H_2(P;\Z)$
on which any closed connection form is nonintegral. But if every 
$M$-bundle with 
closed integral connection form is pulled back from some universal 
bundle ${\mathcal P}\to \Bb$ with this property, then $P\to V$ would  
also be such a pullback, and hence would also have a closed integral 
connection form.  
Thus, if $\Symp(M,\om)$ contains  elements $h_1,h_2$ as above  {\it 
there is no 
universal $M$-bundle ${\mathcal P}\to \Bb$ with closed integral 
connection form.}

A similar argument applies whenever there are two splittings $s_i, 
i=1,2,$
such that the images of $\Ham^{s_i\Z}$ in $\pi_0(\Symp)$ are 
different.  
It follows from Lemma~\ref{le:tau} that this happens iff there 
are integral lifts $\tau_i, i = 1,2$, of $[\om]$ that are stabilized 
by different  subgroups of $\pi_0(\Symp)$.
 To get around this difficulty, one must reformulate
 the classification 
 problem: see Gal--K\c edra~\cite{GalK}.  
\end{example}

We next explain
    a very natural way to think of  a splitting $s^{\Z}$  of $\pi_\Z$
    when $\Poz=\Z$.  Denote by $\tau\in H^2(M;\Z)$ some integer lift 
    of $[\om]$ and by $\rho: L_{\tau}\to M$ the corresponding
    prequantum complex  line bundle.  
Choose a connection $1$-form $\al$ on $L: = L_{\tau}$ with curvature 
$d\al=-\rho^*(\om)$.
Then $\al$ determines a splitting 
$$
s: = s_{\al}
$$
as follows.  By
Stokes' theorem the $\al$-holonomy round a nullhomologous
loop $\ell$ in $M$ 
is multiplication by 
\begin{equation}\labell{eq:almono}
\exp(2\pi i\int_W\om),\quad\mbox{ where }\;\p W = \ell.
\end{equation}
(To see this observe that any such $W$ can be cut open until it is a 
disc and so can be lifted to a $2$-disc $\Tilde W$ in $L$ whose 
boundary projects to the union of $\ell$ with some arcs that are 
covered twice, once in each direction.)
Therefore the holonomy of $\al$ defines a homomorphism 
$m_{\al}$ from the group $Z_1(M)$  of integral 
$1$-cycles in $M$ to $\R/\Z$ that factors through  $SH_1(M,\om;\Z)$.
Now set $s_\al[\ell]$ to be the unique element in $\pi_\Z^{-1}[\ell]$ 
in the kernel of $m_{\al}$.  In other words, we choose $s_\al[\ell]$
so that the $\al$-holonomy round this loop vanishes.
Two connection forms $\al, \al'$ on $L_\tau$ differ by the pullback 
of a closed $1$-form $\be$ on $M$.  Hence if $[\ell]$ is a torsion 
class 
the element $s_\al[\ell]$ is independent of the choice of $\al$.
Therefore each integer lift $\tau$ of $[\om]$ determines a family of 
splittings $s: = s_\al$ that give rise to the same element $C_s\in 
\Hom(\Tor,\R/\Z)$.

\begin{defn}\labell{def:scanon} We say that a splitting 
 $s: = s^{\Z}$ is {\bf $\tau$-canonical} if it has the form $s_{\al}$ 
for some connection form $\al$ on the prequantum bundle 
$\rho:L_{\tau}\to M$.
\end{defn}

\begin{lemma}\labell{le:tau} {\rm (i)} Each splitting $s$ is 
$\tau$-canonical  
    for a unique bundle $L_{\tau}\to M$.
    \sss
    
    \NI
    {\rm (ii)}  Two connections forms $\al, \al'$ on $L\to M$ give 
    rise to the same splitting 
 iff $\al - \al'$ is the pullback 
of an exact $1$-form on $M$.\sss

\NI
{\rm (iii)}  Any two $\tau$-canonical splittings 
    $s,s'$ of $\pi_{\Z}$ are homotopic.
Moreover the corresponding groups 
 $\Ham^{s\Z}$ and $\Ham^{s'\Z}$  are conjugate  by an element of 
 $\Symp_0$.  \sss

\NI
{\rm (iv)}
If $s$ is $\tau$-canonical, the image of $\pi_0(\Ham^{s\Z})$ in 
$\pi_0(\Symp)$  is the stabilizer of $\tau$ in $\pi_0(\Symp)$.
\end{lemma}
\begin{proof} 
Statement (ii) is immediate from the construction.  
To prove (iii) note that two connection forms $\al, \al'$ on $L_\tau$ 
differ by the pullback of a closed $1$-form $\be$ on $M$.  Hence 
the corresponding splittings $s_\al$ and $s_{\al'}$ can be joined by 
a path of splittings that are constant on the torsion loops.  The 
proof of
Lemma~\ref{le:split1} shows that this path 
can be lifted to an isotopy in $\Symp_0$. 
Therefore  $\Ham^{s\Z}$ and $\Ham^{s'\Z}$ are conjugate as in 
Lemma~\ref{le:conj}.

To prove (i) note that the 
set of integer lifts $\tau$ of $[\om]$ is a coset of
 the torsion subgroup $\Tor H^2(M;\Z)$, while the set of splittings 
of $\pi_\Z$ 
 is a coset of
$\Hom (H_1(M;\Z);\R/\Z)=H^1(M;\R/\Z)$.  By (iii) we have set up a 
correspondence $\tau\mapsto s_\al$ between the set of integer lifts 
and the components of
the space of splittings.  Since these are finite sets with the same 
number of elements, we simply have to check that this correspondence  
is injective.  In other words, we need to see that the isomorphism 
class of $L$ is determined by the set of loops that are  
homologically torsion and have
 trivial $\al$-holonomy (where $\al$ is any connection $1$-form.)
But this is an elementary fact about complex line bundles.  In fact, 
given bundles $L, L'$ with connections $\al,\al'$ that have the same 
curvature
and have trivial holonomy round a set of loops generating 
$\Tor(H_1(M;\Z))$, one can adjust $\al'$ so that the monodromies 
agree on a full set of generators for $H_1(M;\Z)$ and then construct 
an isomorphism between the two bundles  
by parallel translation.

Finally note that by Lemma~\ref{le:hsz0}
the image of $\pi_0(\Ham^{s\Z})$ in $\pi_0(\Symp)$ is the kernel of 
$C_s$. 
The proof of (i) shows that if $s$ is $\tau$-canonical then $g\in 
\ker C_s$ iff 
$g^*(\tau) = \tau$.  This proves (iv). 
\end{proof}

We end this section by showing that  when $s$ is $\tau$-canonical the 
group $\Ham^{s\Z}$ has a 
natural geometric interpretation in terms of the bundle $L_{\tau}\to 
M$. 
As in Gal--K\c edra~\cite{GalK},
consider the group $\Gg: = \Gg_{\tau}$  of  all $S^1$-equivariant 
automorphisms of the prequantum line bundle
$L: = L_{\tau}$ that cover a symplectomorphism of  $M$.
Then there  is a fibration sequence of groups and group homomorphisms
$$
{\rm Map}(M,S^1)\to \Gg_{\tau}\stackrel{\rho}\to \Symp_\tau,
$$
where the elements of ${\rm Map}(M,S^1)$ act by rotations on the fibers and 
$\Symp_\tau$  consists of all elements in $\Symp$ that fix the 
given lift
$\tau\in H_2(M;\Z)$ of $[\om]$.  Thus $\Symp_\tau$ is
 a union of components of $\Symp$.  (Note that every $g\in \Symp_\tau$
 does lift to an element in $\Gg_\tau$.  Indeed, since $g^*(\tau) = \tau$ there is a bundle isomorphism
 $\psi:L\to g^*(L)$.  But this induces an isomorphism
 $\psi(x): L_x\to (g^*L)_x = L_{gx}$ for all $x$ and so is a lift of $g$.)
 Gal--K\c edra point out that the line bundle $L_\tau\to M$
extends to a line bundle over
 the universal $M$-bundle over $B\Gg_{\tau}$. Hence
 any symplectic $M$-bundle whose structural group lifts to 
 $\Gg_\tau$ has an integral connection form.
 
Fix a unitary connection $\al$ on $L: = L_\tau$ and
consider the subgroup $\Gg_{\al}\subset \Gg_\tau$  of all $S^1$-equivariant 
automorphisms of the prequantum line bundle
$L$ that preserve $\al$.  If $s$ is the splitting defined by $\al$ then it follows immediately from the definitions  that $\Gg_{\al}$ projects onto
the subgroup $\Ham^{s\Z}$: cf. Kostant~\cite[Theorem~1.13.1]{Ko}.  Hence
there is a commutative diagram 
$$
\begin{array}{ccccc}
S^1&\to& \Gg_{\al}&{\to}& \Ham^{s\Z}\\
\downarrow&&\downarrow&&\downarrow\\
{\rm Map}(M,S^1)&\to&  \Gg_{\tau}&\stackrel{\rho}\to&  \Symp_\tau,
\end{array}
$$
in which the top row is a central extension and the kernel $S^1$ 
of this extension is included in the bottom row as the constant maps.

Now consider the quotient group 
$
\ov\Gg_{\tau}: =\Gg_{\tau}/\Nn
$
 where $\Nn\subset 
{\rm Map}(M,S^1)$  consists
of all null homotopic maps $M\to S^1$.  
The fibration over $\Symp_\tau$ descends to give an extension
\begin{equation}\labell{eq:cext} 
    H^1(M;\Z) \to \ov\Gg_\tau \to \Symp_\tau.
\end{equation}
We give a formula for its defining $2$-cocycle in
Remark~\ref{rmk:eps} below. 
%The next result explains its relation to
%$\Ham^{s\Z}$.  
%Observe also that the subgroup $\Gg_\al$ of $\Gg_\tau$ consisting of automorphisms that fix a given connection form $\al$ intersects the kernel $\Map(M,S^1)$ of the projection $\Gg_\tau\to \ov\Gg_\tau$
%in a copy of $S^1$.  Thus there is a central 
%extension $S^1\to \Gg_\al\to
%\ov\Gg_\tau$, whose restriction to the identity component 
%was considered by Kostant~\cite{Ko}.  (See also Haller--Vizman~\cite[\S3]{HV}.)

\begin{prop}\labell{prop:Z}  Suppose that $\Poz = \Z$.
Given any $\tau$-canonical splitting $s$, 
the inclusion $\Ham^{s\Z}\to\Symp_\tau$ lifts to
a continuous group homomorphism 
$$
\io: \Ham^{s\Z}\to\ov\Gg_\tau
$$
that is a homotopy equivalence.
\end{prop}

\begin{proof}
Since the induced map $\Gg_\al\to \ov\Gg_\tau$ factors through the quotient $\Gg_\al/S^1\cong \Ham^{s\Z}$, we have a commutative diagram
$$
\begin{array}{ccccccc} 
&& \Ham^{s\Z}&\hookrightarrow& \Symp_\tau&\stackrel{\wF_s}\to & BA\\
&&\io\downarrow&&\downarrow&&\\
A &\to &\ov\Gg_\tau &\to& \Symp_\tau&&
\end{array}
$$
where $A: = H^1(M;\Z)$.  Since $\io$ induces an isomorphism on $\pi_0$ by construction, it suffices to check that the boundary map $\p: \pi_1(\Symp_\tau) \to
\pi_0(A) = H^1(M;\Z)$ in the long exact sequence for the bottom row
is given by the usual flux 
$
\Flux: \pi_1(\Symp_0)\to \Ga\subset H^1(M;\Z). 
$
This is well known, but we include the proof for completeness.

%\begin{proof}  Choose a connection $\al$ on $L$.
The space $\Ppp_*\Symp_0$ of 
based paths $(h_t)$ in $\Symp_0$ acts on $L$ by
 taking $(x,\theta)$ to its image under $\al$-parallel translation 
 along the path $h_t(x)$.  If $(h_t) $ and $(h_t')$ are two  
 paths with endpoint $h$, then
 $$
(h_t)\cdot(x,\theta) =  \la_f\bigl((h_t')\cdot(x,\theta)\bigr)=
(h(x),\theta+f(x)),
 $$
 where  $f(x)$ is the 
 holonomy of $\al$ around the loop  based at $h(x)$
 that is formed by first going back along
 $h_t'(x)$ and then forward along $h_t(x)$.  When $(h_t) $ and 
$(h_t')$ 
 are homotopic,  the function $f:M\to S^1 = \R/\Z$ is obviously 
 null homotopic.  
Therefore this gives a homomorphism 
%%%
%\begin{equation}\labell{eq:tj}
$$
\Tilde\jmath: \Tilde{\Symp}\to 
\ov\Gg_\tau.
$$
%\end{equation}
%%%

Consider the loop
$\phi=\{h_t\circ (h_t')^{-1}\}$, and observe that 
$f:M\to S^1$ is homotopic to the function $f_\phi$ defined by
  $$
 f_{\phi}(x) = \mbox{$\al$-holonomy round } \phi_t(x),\qquad x\in M.
 $$
 Therefore $f^*[ds] = f_\phi^*[ds]$, and so it suffices to prove
that for any loop  $\ell$ in $M$ 
%%%
%\begin{equation}\labell{eq:fphi}
$$
 \int_{\ell} f_{\phi}^*(ds)  = \Flux(\phi)[\ell].
 $$
%\end{equation}
%%%
But  if $\Tilde f_\ell:[0,1]\to \R$ is any lift of the function 
$
S^1\to S^1:\; s\mapsto f_{\phi}(\ell(s)),
$
 equation (\ref{eq:almono}) implies that
  $$
  \Tilde f_\ell(1)-\Tilde 
f_\ell(0) = \int_{[0,1]\times S^1}\Psi^*\om=\Flux(\phi)[\ell],
$$
where $\Psi(t,s) = \phi_t\ell(s)$. The result follows.
\end{proof}

Gal--K\c edra give a different proof of the above result in \cite[\S4]{GalK}.

%%%%%%%%%%%%%%%%%%%%%%%%%%%%%%%%%%%%%%%%%%%%%%%%%%%%%%%%%%%%%%%%%%%%%%%%
\subsection{Lifting $\wF_s$}\labell{sec:lift}
%%%%%%%%%%%%%%%%%%%%%%%%%%%%%%%%%%%%%%%%%%%%%%%%%%%%%%%%%%%%%%%%%%%%%%%%

We now investigate situations in which the flux homomorphism 
extends
to a continuous crossed homomorphism $\Tilde F$ with values in 
$H^1(M;\R)/\Ga$.
One obvious situation when this happens is when $\Ga = H^1(M;\Poz)$
and $H_1(M;\Z)$ has no torsion
since then we can take $\TF = \wF_s^\Z$.  For example, we can take 
$M$ to be
the product of the standard $2n$-torus with a simply connected 
manifold.
The other cases in which we know that
$\TF$ exists have $\Ga=0$,
so that in fact $\wF_s$ lifts to
 $H^1(M;\R)$.

As observed by Kotschick--Morita~\cite{KM}, 
such a lift exists in the  monotone case.  In this case $\om$ always 
has a closed extension --- by a multiple of the first Chern class of 
the vertical 
tangent bundle --- and so the flux group vanishes by 
Proposition~\ref{prop:noact}.
We now explain their argument.

Kotschick and Morita consider the symplectomorphism group $\Symp^\de$ 
with the  discrete topology and observe that the universal 
$M$-bundle $M_{\Symp^{\de}}\to B\Symp^\de$ has a natural flat 
connection (or foliation) that 
is transverse to the fibers.  
Therefore the fiberwise symplectic 
form $\om$ has a closed extension $\Tilde\om$ given by the associated
 connection form. Suppose now that $(M,\om)$ is monotone, i.e. that 
the 
symplectic class $[\om]\in H^2(M)$ is  a nonzero multiple $\la 
c_1(TM)$ 
of the
first Chern class of the tangent bundle $TM$.
Then $[\om]$  has another extension $\la v$
 where $v: = c_1^{\rm Vert}$ denotes the first Chern class of the 
vertical 
 tangent bundle.  By construction the class $[\Tilde\om] -\la v\in 
 H^2(M_{\Symp^{\de}})$ vanishes when restricted to the fiber.  
 Therefore it projects to a class 
 $$
 [F_{KM}]\in E_{\infty}^{1,1} =
E_2^{1,1} =  H^1(B\Symp^\de; H^1(M;\R)) = H^1_{\EM}(\Symp^\de; 
H^1(M;\R)),
 $$
 where $E_{\infty}^{1,1}$ is the $(1,1)$ term in the Leray--Serre 
 spectral sequence.   
 Each class $[F]\in H^1(B\Symp^\de; \Aa)$ may be represented by a  
crossed 
 homomorphism $F:\Symp^\de\to \Aa$ whose restriction to the subgroup
 $\Symp_H^\de$ that acts trivially on $\Aa$ is unique.  Therefore it 
 makes sense to talk about the restriction of  $F_{KM}$ to the 
 identity component $\Symp_0^{\de}$.
 Kotschick and Morita show that this is just the usual flux 
 homomorphism.\footnote
 {
 As pointed out by K\c edra--Kotschick--Morita~\cite{KKM} this argument works whenever  $[\om]$ extends to $H^2(M_{\Symp})$, which leads to an alternative proof of Corollary~\ref{cor:obs}. See also 
 Remark~\ref{rmk:KKM}.
 }

 We now show that the Kotschick--Morita class
 $[F_{KM}]$ also lifts the class $[\wF_s]$.  We do this by giving an explicit formula for a representative $\Tilde F_s$ of $[F_{KM}]$.
 The task here is to use the monotonicity of $(M,\om)$ to define a correction term for $\wF_s$ so that it takes values in $H^1(M;\R)$ rather than in $H^1(M;\R/\Poz)$.
 If $M$ is a Riemann surface, this is quite easy: 
 see Remark~\ref{rmk:exten}(i).  However, one has to work a bit harder to derive the general formula.   We shall use the 
 notation of \S\ref{ss:Z}.
 
 \begin{thm}\labell{thm:TF}  Suppose that $(M,\om)$ is monotone, 
i.e. that the 
symplectic class $[\om]\in H^2(M)$ is  a nonzero multiple $\la c_1$ 
of the
first Chern class $c_1: = c_1(TM)$.  Choose a splitting $s$ that is 
$\tau$-compatible where $\tau = [c_1]\in H^2(M;\Z)$.
Then  $\widehat F_s:\Symp(M,\om)\to \Aa$
lifts to a crossed homomorphism 
$$
\Tilde F_s: \Symp(M,\om)\to H^1(M;\R).
$$
Moreover this lift represents 
$[F_{KM}]$ in $H^1(B\Symp^\de; H^1(M;\R))$.
\end{thm}

\begin{proof}  Fix a basepoint $x_0$ in $M$ and denote by
$\Symp_*$ (resp. $\Ham_*$) the subgroup of $\Symp$ 
(resp. $\Ham$) that fixes $x_0$.  We  first define $\Tilde F_s$
on the subgroup $\Symp_*$ by looking at the action of $\Symp_*$ on a certain line bundle over the universal cover of $M$.  We then show that $\Tilde F_s$ lifts 
$\widehat F_s$, and  that  $\Tilde F_s$ vanishes on 
$\Ham_*$.  It is then easy to extend $\Tilde F_s$ to a crossed homomorphism on $\Symp$ that lifts $\widehat F_s$: see Step 4.
 The proof that $\Tilde F_s$ represents 
the class $F_{KM}$ is based on giving an adequate description for the 
vertical first Chern class $v$. \MS

\NI
{\bf Step 1:}  {\it Definition of $\Tilde F_s$ on $\Symp_*$.}
\smallskip

By rescaling $[\om]$ we may suppose that $\Poz = \Z$. 
    Then $c_1(TM) = N[\om]$ where $|N|$ is the minimal  Chern number.
Denote by $\Tilde M\to M$ the cover of $M$ associated to the 
homomorphism
$\pi_1(M)\to H_1(M;\Z)$, let $L\to M$ be the complex line bundle 
with first Chern class $c_1(TM)$ and denote by $\Tilde L\to \Tilde 
M$  its 
pullback to $\Tilde M$.     Denote by
$\Tilde x_0 = (x_0,[\ga_0])$ the  base point in 
$\Tilde M$ corresponding to $x_0$,
where $\ga_0$ denotes the constant 
path at $x_0$.

Fix a (unitary) connection $\al$ on $L$, and let
$\Tilde\al$ be its pullback to  $\Tilde L$. (Thus $d\Tilde\al$ is the pullback of $-N\om$.)
Denote by $\Gg_{\Tilde \tau}$ the group of 
$S^1$-equivariant diffeomorphisms of $\Tilde \tau$ that cover 
symplectomorphisms of $M$. 

 For each $g\in \Symp_*$ define 
$
\Tilde g: \Tilde L\to \Tilde L
$
as follows:  the points of $\Tilde M$ are pairs $(x,[\ga_x])$ 
where $[\ga_x]$ is an equivalence class of paths from $x_0$ to $x$
and we set
$$
\Tilde g(x,[\ga_x],\theta): = (gx, [g\ga_{x}], \theta')
$$  
where 
$\theta'\in \Tilde L_{(gx, [g\ga_{x}])}$ is the image of
$\theta\in \Tilde L_{(x, [\ga_{x}])}$ under 
$\Tilde\al$-parallel translation, first back along $\ga_x$ to $x_0$ 
and then 
forwards along $g(\ga_x)$.  This is independent of the choice of 
representative $\ga_x$ for $[\ga_x]$ because $g\in \Symp_*$.
(Use equation~(\ref{eq:almono}) and the fact that $g$ preserves
$\om$.) 
Further the map $g\mapsto \Tilde g$ is a group homomorphism.
%(The bundle $\Tilde L\to \Tilde M$ is independent of the
%choice of $\tau$  because $H_1(\TM) = 0$.) 

By construction  $L$ is isomorphic to $\La^n(T^*(M))$, where $T^*M$ is 
given a complex structure compatible with $\om$.  Hence each $g\in \Symp$ lifts to a bundle automorphism $g_L$ of $L$ that can be chosen to preserve the Hermitian structure of $L$.  Therefore the $1$-form $\Tilde g^*(\Tilde\al)- \Tilde\al$ is the pullback of the $1$-form $(g_L)^*(\al)-\al$ on $M$ and hence is exact. 
For $g\in \Symp_*$  we define $f_g:\Tilde M\to \R$  
to be the unique  function such that $f_g(x_0)=0$ and 
$$
\Tilde g^*(\Tilde\al) = \Tilde\al - df_g.
$$
Then set 
$$
\Tilde F_s(g)([\ell]):= \frac 1N \Bigl(f_g(\Tilde\ell(1))-
f_g(\Tilde\ell(0))\Bigr)\in \R,\qquad g\in \Symp_*,\;[\ell]\in 
H_1(M;\Z)/\Tor,
$$
where $\ell$ is a based loop in $M$ representing $[\ell]$
and the path $\Tilde\ell$ is its lift to $\Tilde M$ with initial point $\Tilde x_0$.
In other words,
$$
\Tilde F_s(g)([\ell]):= \frac 1N \Bigl(m_{\Tilde g^*\Tilde\al}(\Tilde\ell) - 
m_{\Tilde\al}(\Tilde\ell)\Bigr),
$$
where $m_{\Tilde\be}(\Tilde\ell)$ denotes the
holonomy  of the connection $\Tilde\be$ along the path 
$\Tilde\ell$.  
(This is only defined mod $\Z$, but the difference between two
connections gives an element in $\R$.)  It is easy to check that 
$\Tilde F_s(g)([\ell])$ does not depend on the chosen representative
$\ell$ for $[\ell]$.  This completes Step 1.\QED

\NI
{\bf Step 2:} {\it $\Tilde F_s$ reduces  mod $\Z$ to 
$\wF_s$.}\smallskip

It is immediate from the definitions that 
$
\Tilde F_s:\Symp_*\to \R$ is a crossed 
homomorphism. We claim that its mod $\Z$ reduction coincides with
$\wF_s$.  Suppose first that $g[\ell] = [\ell]\in H_1(M;\Z)$.  Then the
lifted paths $\Tilde\ell$ and $\Tilde g(\Tilde\ell)$ have the same
endpoint and form the boundary of a $2$-chain $\Tilde C$ in $\Tilde M$
that lifts a chain $C$ in $M$ with boundary $C = g(\ell) - \ell$. 
Hence we can calculate the mod $\Z$-reduction of $\Tilde
F_s(g)([\ell])$ as follows:
\begin{eqnarray*}
\Tilde F_s(\Tilde g)([\ell]) & = & \frac 1N\Bigl(m_{\Tilde 
g^*\Tilde\al}(\Tilde\ell) - m_{\Tilde\al}(\Tilde\ell)\Bigr)\\
& = & \frac 1N\Bigl(m_{\Tilde\al}(\Tilde g(\Tilde\ell)) -
m_{\Tilde\al}(\Tilde\ell)\Bigr)\\
& = & \frac 1N \int_{\p\Tilde C} -\Tilde \al\; =\;
 \int_C\om \\
& = & \wF_s(g)[\ell] \in \R/\Z.
\end{eqnarray*}
In general, $\wF_s(g)[\ell]$ is defined to 
be the area of a chain with boundary  $g(\ell) - \ell'$, where $\ell, 
\ell'$ are chosen to have zero $\al$-holonomy (i.e. so that they 
project to elements in ${\rm Im\,}s$) and $\ell'$ is homologous to 
$g(\ell)$.   But now $\Tilde g(\Tilde\ell)$ and $\Tilde\ell\,\!'$ 
have the same 
endpoint and so
\begin{eqnarray*}
\Tilde F_s(\Tilde g)([\ell]) & = & \frac 1N\Bigl(m_{\Tilde 
g^*\Tilde\al}(\Tilde\ell) - m_{\Tilde\al}(\Tilde\ell\,\!')\Bigr),\\
& = & \frac 1N\Bigl(m_{\Tilde\al}(\Tilde g(\Tilde\ell)) -
m_{\Tilde\al}(\Tilde\ell\,\!')\Bigr),\\
& = & \wF_s(g)[\ell] \in \R/\Z,
\end{eqnarray*}
as before.
Thus  $\Tilde F_s$ lifts $\wF_s$  on $\Symp_*$. \QED

\NI
{\bf Step 3:} {\it $\Tilde F_s$ vanishes on $\Ham_*$.}\smallskip

Note first that although
$\Ham$ is a perfect group, its subgroup $\Ham_*$ is not --- it
supports the Calabi type homomorphism $h\mapsto \int_0^1\int_M
H_t\om^n$, where $H_t$ generates $h$ and is normalized by the
condition that $H_t(x_0) = 0$.  Hence this step does require proof.

To this end, let $h\in \Ham_*$ and choose any path $h_t$ 
from the identity to $h: = h_1$ such that the 
$\al$-holonomy round the loop
$h_t(x_0)$ is trivial.  Then $h_t$ lifts to a path in $\Tilde M$ and
thence  to a
path $\Hat h_t: \Tilde L\to \Tilde L$
given   by taking the 
$\Tilde\al$-holonomy along $h_t(x)$.  Let $s\mapsto \ell(s)$ be  a 
based loop in $M$ with lift $\Tilde\ell$. 
Consider the map 
$$
Y_{\eta}:I^2\to \Tilde M,\quad (s,t)\mapsto \Hat h_t(\Tilde\ell(s)),
$$
and trivialize the bundle $Y^*\Tilde L$ by parallel translation along 
the horizontal line $t=0$ and the verticals $s=const.$  
Then in this trivialization the 
map 
$$
\Hat h_1:\Tilde L_{Y([0,1]\times \{0\})}\to L_{Y([0,1]\times \{1\})}
$$ 
has the form
$(s,0,\theta)\to (s,1,\theta)\in Y^*(\Tilde L)$
while the corresponding map defined by
$\Tilde h: = \Tilde{\io}(h)$  has the form
$$
\Tilde h: (s,0,\theta) \mapsto (s,1,\theta- f(s))
$$
where
$f(s)$ is the area of the rectangle $Y([0,s]\times [0,1])$.
(Here we have used the fact that  the $\al$-holonomy round the loop
$h_t(x_0)$ is trivial so that our trivialization gives the obvious 
identification of the fiber $L_{x_0}$ over $(0,0)$ with that (also 
$L_{x_0}$) over $(0,1)$.  Further, because $s\mapsto
Y(s,0,\theta)$ is $\Tilde\al$-parallel,
the path $s\mapsto \Hat h(Y(s,0,\theta))$ 
is also $\Tilde\al$-parallel 
by definition of $\Hat h$.
  Therefore, in this trivialization the restriction of
$\Tilde h^*\Tilde\al$ to $\Tilde\ell$ 
is $\Tilde\al - df(s)$. Hence $\Tilde F_s(h)[\ell] = f(1)$. 
But because  $h_t$ is a Hamiltonian path with no flux through $\ell$ 
we must have $f(1) = 0$.  \QED

\NI
{\bf Step 4:} {\it End of the proof of the first statement.}\smallskip

 Given $g\in \Symp$ we now set
$$
\Tilde F_s(g): = \Tilde F_s(hg),
$$
where $h\in \Ham$ is any element such that $hg\in \Symp_*$. 
This is independent of the choice of $h$ because $\Tilde F_s=0$ on
$\Ham_*$.  Moreover $\Tilde F_s:\Symp\to H^1(M;\R)$ is a
crossed homomorphism.  It lifts $\wF_s$, since this also vanishes on
$\Ham$.  This completes the proof of the first statement in the proposition.\QED

\NI
{\bf Step 5:} {\it  $\Tilde F_s$ represents the class $[F_{KM}]$
in $H^1_{\EM}(\Symp^\de;\R)$.}\smallskip

 This is easy to see on the subgroup 
$\Ham$ since both crossed homomorphisms vanish there. Hence it 
suffices to check this statement for $\Symp_*$.
The difficulty in proving this is to find a suitable way to
calculate the class $v: = c_1^{\rm Vert}$.

   We first consider 
the subgroup
$S_{*H}: = \Symp_{*H}$  of $\Symp_*$ that acts trivially on 
$H_1(M;\Z)$. 
The $E^2_{1,1}$-term of the integral homology 
spectral sequence for the pullback of 
$M_{\Symp^{\de}}\to 
B\Symp^{\de}$ to $BS_{*H}^{\de}$
is isomorphic to the product $H_1(B S_{*H}^{\de})\otimes 
H_1(M;\Z)$.  It  is generated by cycles 
$$
C_{g,\ell}: = Z_{g,\ell}\cup W_{g,\ell},
$$
where $W_{g,\ell}$ is a $2$-chain in 
the fiber $M_*$ over the base point with boundary $g(\ell)-\ell$,
and $Z_{g,\ell}$  is a cylinder  
lying over the loop in the base corresponding to $g\in 
\Symp^{\de}$  with boundary $\ell - g(\ell)$.  Since the class $[\tom]$ 
vanishes on $Z_{g,\ell}$,
$$
\int_{C_{g,\ell}}\tom = \int_{W_{g,\ell}} \om.
$$
(Note that the mod $\Z$ reduction of this class is $\wF_s(g)[\ell]$ 
as 
one would hope.)
Since $S_{*H}$ does not act on $L$ but does act on the pullback 
$\Tilde L$, 
to understand the vertical Chern class $v = c_1^{\rm Vert}$ 
we should think of $M$ as the 
quotient of $\Tilde M$ by the group ${\rm G}: = H_1(M;\Z)$ and
consider the 
corresponding ${\rm G}$-equivariant pullback 
line bundle 
$
\Tilde L_{S^{\de}}\to \Tilde M_{S^{\de}}.
$ 
(Here we denote $S^\de : = S^\de_{*H}$.)
Then $C_{g,\ell}$ lifts to a cycle  $\Tilde C_{g,\ell}: = \Tilde 
Z_{g,\ell}\cup 
\Tilde W_{g,\ell}$ in $\Tilde M_{S^{\de}}$.  There is a 
trivialization 
of $\Tilde L_{S^{\de}}$ over $\Tilde Z_{g,\ell}$ that restricts  to 
$\Tilde\al$ over
$\Tilde{\ell}$ and to $(\Tilde g)_*(\Tilde\al)$ over $\Tilde 
g(\Tilde{\ell})$. 
Extend this to any connection $\Tilde\be$  of $\Tilde L$ over $\Tilde 
W_{g,\ell}$.
Then $v(C_{g,\ell})$ is given by integrating the curvature of 
this connection over $\Tilde C_{g,\ell}$.  Since this connection is 
flat 
over $\Tilde Z_{g,\ell}$, the relevant part of the integral is over 
$\Tilde W_{g,\ell}$.  Then because $-N\tom$ pulls back to $d\Tilde\al$
\begin{eqnarray*}
\Bigl\langle \tom - \frac 1N v, C_{g,\ell}\Bigr\rangle & = &
\frac 1N\int_{\Tilde W_{g,\ell}} (-d\Tilde\al + d\Tilde\be)\\
&=& 
\frac 1N \Bigl(\int_{\Tilde g(\Tilde\ell)} \Tilde\be -\Tilde\al - 
\int_{\Tilde\ell} \Tilde\be -\Tilde\al
\Bigr) \\
& = & \frac 1N  \int_{\Tilde g(\Tilde\ell)} \Tilde\be -\Tilde\al\;\; =\;\;
\frac 1N \int_{\Tilde\ell} (\Tilde g)^*(\Tilde\be -\Tilde\al)\\
&=& 
\frac 1N \int_{\Tilde\ell} \Tilde\al -(\Tilde g)^*\Tilde\al
\;\;=\;\; \Tilde F_s (g)[\ell].
\end{eqnarray*}

To extend this argument to the full group $S: = \Symp_*$ 
we consider the $2$-{\it chains} 
$$
{C'_{g,\ell}} = Z_{g,\ell}\cup W'_{g,\ell},
$$
where $\ell$ is now a loop with trivial $\al$-holonomy,
$Z_{g,\ell}$ is as before, and $W'_{g,\ell}$ is any $2$-chain in 
$M_*$ 
with boundary $g(\ell) - \ell'$, where $\ell'$ is homologous to 
$g(\ell)$ 
and also has zero $\al$-holonomy. Thus $\p C'_{g,\ell} = \ell -
\ell'$,
where both loops have zero $\al$-holonomy.
 
Any element in $E^{1,1}_2 = H^1(BS^{\de}; H^1(M;\R))$ is determined 
by 
its values on the integral $2$-cycles $Z$ that are sums of chains 
of the form
$C'_{g,\ell}$ with chains in the fiber $M_*$ whose boundary consists 
of sums of
 loops  with trivial $\al$-holonomy. 
Again we lift each such cycle to a cycle $\Tilde Z$ in 
$\Tilde M_{S^\de}$, and evaluate
 $v$ on $\Tilde Z$ by integrating the curvature of a
 suitable connection form for the pullback of $\Tilde L$ to $\Tilde 
Z$.
 As before we suppose that this connection equals $\Tilde\al$ on the
 \lq\lq free'' boundary arcs $\Tilde\ell, \Tilde\ell'$.  Since these 
have
 trivial $\Tilde\al$-holonomy by construction, the pieces of $\Tilde 
Z$ formed 
 from chains in $M_*$ do not contribute to $[\om]-\frac 1N v$, while 
  the contribution of $C'_{g,\ell}$ is 
$$
\Bigl\langle \tom - \frac 1N v,\, C'_{g,\ell}\Bigr\rangle
= \Tilde F_s (g)[\ell]
$$
as before.   This completes the proof.
\end{proof}

\begin{rmk}\labell{rmk:exten}\rm  (i)  If $(M,\om)$ is a Riemann surface, then there is another simpler description
of $\TF_s$. In the notation of \S2, $\TF_s(h)(\la) = \int_C\om - \frac 1N c_1(TM|_C)$ where $C$ is any integral $2$-chain whose boundary 
represents $h_*(s\la) - s(h_*\la)$ and  $c_1(TM|_C)$ is the relative Chern number of the restriction of $TM$ to $C$ with respect to the obvious trivializations of $TM$ along the (embedded) boundary of $C$.
\MS

\NI
(ii)  The above argument used in an essential way the fact that every symplectomorphism lifts to 
an automorphism of $L = \La^n(T^*M)$.  The construction of $\TF_s$ works whenever there is a line bundle $L'$ over the universal $M$-bundle 
$M_{\Symp}\to B\Symp$ such that $c_1(L')$ restricts on the fiber to an integral lift $\tau$ of $[\om]$, i.e. whenever there is a lift such that $\Symp_\tau = \Symp$.  The second part of the argument, showing that $\TF_s$ represents the corresponding element in the $E^{1,1}_2$-term of the spectral sequence also goes through.
Hence, under these circumstances, our construction gives an explicit
formula for the crossed homomorphism whose existence is established by
K\c edra--Kotschick--Morita in~\cite[Theorem~6]{KKM}.
\end{rmk}

%%%%%%%%%%%%%%%%%%%%%%%%%%%%%%%%%%%%%%%%%%%%%%%%%%%%%%%%%%%%
\subsection{The atoroidal case}\labell{ss:ator}
%%%%%%%%%%%%%%%%%%%%%%%%%%%%%%%%%%%%%%%%%%%%%%%%%%%%%%%%%%%%

We now consider the situation when $\om$ vanishes on tori and/or spheres.   Recall from \S1 that the group $\Symp_\pi$ consists of all symplectomorphisms that are isotopic to an element in $\Symp_*$ that acts trivially on $\pi_1(M,x_0)$, where $x_0$ is the base point in $M$.  This is equivalent to saying that the elements in $ \Symp_\pi\cap \Symp_*$ act on 
$\pi_1(M,x_0)$ by inner automorphisms.  We denote by 
$$
\Symp_{*\pi}
$$
 the subgroup of $\Symp_\pi\cap \Symp_*$ that acts trivially on $\pi_1(M,x_0)$.  The (disconnected) group $\Ham_*$ is a subgroup of $\Symp_{*\pi}$ because the evaluation map $\pi_1(\Ham)\to \pi_1(M,x_0)$ is trivial: see~\cite{Mcox} for example.  Note also that although $\Ga$ may not vanish when $\om=0$ on $\pi_2(M)$, the Flux homomorphism vanishes
  on loops in $\Symp_{*\pi}\cap \Symp_0$ since its value 
   is then given by evaluating $[\om]$ on $2$-spheres.  Therefore, in this case  Flux is  well defined as a homomorphism 
   $\Symp_{*\pi}\cap \Symp_0\to H^1(M;\R)$.

\begin{lemma}  If $\om =0$ on spheres then $\Flux$ extends to a homomorphism $F_\pi:\Symp_{*\pi}\to H^1(M;\R)$.
\end{lemma}
\begin{proof}
% 
% \NI
% {\bf Proof of Proposition~\ref{prop:symphyp}.} It is immediate that
% $\Ga=0$ when $\om = 0$ on tori.  Hence we need only  show that Flux %extends to
% a homomorphism $F_\pi:\Symp_{\pi}\to H^1(M;\R)$.
%To this end,  fix a base point $x_0\in M$ and
Denote by $S\pi_1(M)$ the group
formed by equivalence classes $\langle\ga\rangle$ of based loops in
$(M,x_0)$ where two loops $\ga_0,\ga_1$ are equivalent iff they may be
joined by a based homotopy $\ga_t$ of zero symplectic area, i.e.
$$
    \int_C \psi^*\om = 0, \qquad
  \psi:C = S^1\times [0,1]\to M, \;(s,t)\mapsto \ga_t(s).
  $$
Because $\om$ vanishes on spheres, the symplectic area of a based homotopy
  between two homotopic loops $\ga_0,\ga_1$ is independent of the choice
  of homotopy.  As in Lemma~\ref{le:split0} it follows readily that
  there is an exact sequence of groups
\begin{equation}\labell{eq:Shtpy}
  0\to \R\to S\pi_1(M)\stackrel{\pi}\to \pi_1(M)\to \{1\}.
\end{equation}

  We now define a map $s: \pi_1(M)\to S\pi_1(M)$ such that $\pi\circ s = \id$.   Since $\om$ vanishes on spheres it lifts to an exact form $\Tilde\om$ on the universal cover  $\TM$ of $M$.  Choose a base point $\tx_0$ of  $\TM$ that projects to $x_0$ and choose a $1$-form $\be$ on $\TM$ such
  that $d\be =\tom$.   Then define $s: = s_{\be}$ by setting
  $s_{\be}([\ga]) =\langle\ov\ga\rangle$ where the based loop $\ov\ga$
  represents $[\ga]$ and is such that $\be$ integrates to zero over the
  unique lift of $\ov\ga$ to a path in $\TM$ starting at $\tx_0$.  Note
  that $s: = s_\be$ need not be a group homomorphism because $\be$ is not
  invariant under deck transformations.

For each $g\in \Symp_{*\pi}$ define a function
$$
F_{\pi}(g): \pi_1(M)\to \R
$$
by setting
$F_\pi(g)([\ga])$ equal to the symplectic area of a based homotopy 
joining 
 $\ov{g(\ga)} = \ov\ga$ to $g\ov\ga$.  
We claim that the map $F_\pi(g): \pi_1(M)\to \R$ is a homomorphism
for each $g\in \Symp_{*\pi}$ and so defines an element $F_\pi(g)$ in $H^1(M;\R)$.  To
see this, 
 define $A(\ga,\de)$ to be the $\om$-area of any homotopy from
  $\ov\ga *\ov\de$ to $\ov{\ga*\de}$, where $\ga,\de\in \pi_1(M)$ and
  $\ov\ga *\ov\de$ denotes the loop that first goes round $\ov\ga$ and then
  $\ov\de$.
 Now
$F_\pi(g)([\ga*\de])$ is the area of any homotopy from
$\ov{g(\ga*\de)}$ to $g(\ov{\ga*\de})$.  Consider the three part
homotopy that goes first from $\ov{g(\ga*\de)}$ to
$\ov{g\ga}*\ov{g\de}$, then to $g\ov{\ga}*g\ov{\de}$ and finally to
$g(\ov{\ga*\de})$.  
The first of these homotopies has area $-A(g\ga,g\de) = -A(\ga,\de)$
since $g$ acts trivially on $\pi_1$.  Therefore this cancels out the
area of the third homotopy.  The middle homotopy can be chosen to be
the juxtaposition of the homotopies used to define $F_\pi(g)[\ga]$ and
$F_\pi(g)[\de]$ and so has area equal to their sum.  Thus
$F_\pi(g)([\ga*\de]) = F_\pi(g)([\ga]) + F_\pi(g)([\de])$.

One now shows that $F_\pi: \Symp_{*\pi}\to H^1(M;\R)$ is a homomorphism as
in Proposition~\ref{prop:wf}.  \end{proof}

\begin{cor}  If $\om$ vanishes on tori then 
$F_\pi$ extends to a  homomorphism $\Symp_{\pi}\to H^1(M;\R)$.
\end{cor}
\begin{proof} 
Since each $g\in \Symp_\pi$ may be joined 
 to a point $gh_1\in \Symp_{*\pi}$ by a Hamiltonian path $h_t\in \Ham$,
 we extend $F_{\pi}$ by setting
 \begin{equation}\labell{eq:hext} F_{\pi}(g) : = F_{\pi}(gh),\qquad g\in
 \Symp_{\pi}, \;\;gh\in \Symp_{*\pi},\;\;h\in \Ham.
\end{equation}
This is independent of the choice of
$h_t\in \Ham$ because $\om=0$ on tori.
\end{proof}

The above proof was written using geometric language, to
imitate the definition of $\wF_{s}$.  However, when $g\in \Symp_{*\pi}$
it is perhaps more
illuminating to write $F_\pi(g)\in \Map(\pi_1(M),\R)$ in the form
\begin{equation}\labell{eq:F}
    F_\pi(g)(\ga) = \ka_g(\Hat\ga),\qquad \ga\in \pi_1(M,x_0),\; g\in
    \Symp_{*\pi}
\end{equation}
where $\ka_g:\TM\to \R$ is the unique function such that
$
\tg^*\be = \be + d\ka_g,$ and $ \ka_g(\tx_0) = 0.
$
Here $\tg$ denotes the lift of $g$ to the universal cover 
$\TM$ that fixes the base point
$\tx_0$, and $\Hat\ga\in \pi^{-1}(x_0)\subset\TM$ is the end point
of the lift of any representative of $\ga\in \pi_1(M)$.  Hence 
$$
\ka_{gh}
= \ka_h +\Th^*\ka_h,\qquad F_\pi(gh) = F_\pi(h) + h^*F_\pi(g).$$

The above identities hold for all $g\in \Symp_*$.
However, when $g\in\Symp_{*\pi}$ we saw above that the map $F(g)$ is a group
homomorphism $\pi_1(M)\to \R$ for each $g$.  In general this is not
true.  To explain the algebra, denote by $\psi_{\ga}:\TM\to\TM$ the
deck transformation corresponding to $\ga\in \pi_1(M)$ and define
the function $f_{\ga}:\TM\to \R$ by
$$
(\psi_\ga)^*\be = \be + df_\ga,\qquad f_{\ga}(\tx_0) = 0.
$$
Then 
\begin{equation}\labell{eq:B}
f_{\ga\de} - f_{\de} - (\psi_{\de})^*f_{\ga} = f_\ga(\Hat\de),
\end{equation}
where we think of the right hand side as a constant function on $M$.
Therefore
\begin{equation}\labell{eq:A}
A(\ga,\de)= \int_{\psi_{\ga}\ov\de}\be = \int_{\ov\de}(\psi_{\ga})^*(\be)
= \int_{\ov\de} df_{\ga} = f_\ga(\Hat\de).
\end{equation}
Hence $A:\pi_1\times \pi_1\to \R$ satisfies the cocycle condition
$$
A(\de,\eps) - A(\ga\de,\eps) + A(\ga,\de\eps) - A(\ga,\de) = 0.
$$
(When checking this,  it is useful to think of $A(\ga,\de)=
(\psi_\eps)^*A(\ga,\de)$ as a constant
function on $\TM$ and to use (\ref{eq:B}).)  Then the identity 
$\Tg\circ \psi_{\ga} = \psi_{g\ga}\circ\Tg$ implies that
$$
f_{g\ga}\circ \Tg - f_{\ga} = \ka_g\circ\psi_{\ga} -\ka_g - \ka_g(\Hat\ga).
$$
It follows that
$$
F_\pi(g)(\ga\de) -F_\pi(g)(\ga) - F_\pi(g)(\de) = f_{g\ga}(\Hat{g\de}) -
f_{\ga}(\Hat\de) = g^*A(\ga,\de) -A(\ga,\de).
$$
Therefore, if we modify (\ref{eq:F}) by setting
$$
F_\pi(g,\ga)(\de):= \ka_g(\Hat\ga) +g^*A(\ga,\de) -A(\ga,\de),
$$
we find that
$$
F_\pi(g,1)(\ga\de) = F_\pi(g,1)(\ga) + F_\pi(g,\ga)(\de).
$$

This discussion applies when $\om$ vanishes on $\pi_2(M)$.  If $\om$
also vanishes on tori\footnote
{
If $\om$ does not vanish on tori, we still have $F(g) = \Flux (g)$ mod
$\Ga$, provided that the evaluation map $\pi_1\Symp_0\to \pi_1(M)$ is
surjective.  So in this case we also get an extension of $\Flux$.}
then $\Ga=0$ and $F_\pi(g) = \Flux (g)$ for $g\in \Symp_{*0}$, so that we can
extend $F_\pi$ to the whole group $\Symp$ by equation (\ref{eq:hext}). 
Hence we have shown:

\begin{prop} Suppose that $\om$ vanishes on tori.  Then Flux
 extends to a crossed homomorphism
$$
F_\pi: \Symp \to \Map(\pi_1(M);\R),
$$
such that $F_\pi(g)\in \Hom(\pi_1(M);\R)$ for $g\in \Symp_{\pi}$.
\end{prop}

Because of the rather complicated algebraic structure of the map
$F_\pi(g)$ when $g$ acts nontrivially on $\pi_1$, the kernel of $F_\pi$
will not in general intersect every component of $\Symp$.  Note also
that if $\om$ lifts to an exact form on the abelian cover of $M$ given
by the homomorphism $\pi_1(M)\to H_1(M;\Z)/\Tor$, one can play the
same game there, defining $F_H$ and $A$ by equations~(\ref{eq:F}) and
(\ref{eq:A}).  Thus:

\begin{cor}\labell{cor:ext} Suppose that $[\om]\in
\bigl(H^1(M;\R)\bigr)\,\!^2\subset H^2(M;\R)$ and that $\om=0$ on tori. 
Then $\Flux$ extends to a crossed homomorphism
$$
F_H:\Symp_H\to H^1(M;\R),
$$
where $\Symp_H$ is the subgroup of $\Symp$ that acts trivially on $H_1(M;\R)$.
\end{cor}

%%%%%%%%%%%%%%%%%%%%%%%%%%%%%%%%%%%%%%%%%%%%%%%%%%%%%%%%%%%%%%%%%%%%%%%%
\subsection{$c$-Hamiltonian bundles and covering 
groups}\labell{sec:cHam}
%%%%%%%%%%%%%%%%%%%%%%%%%%%%%%%%%%%%%%%%%%%%%%%%%%%%%%%%%%%%%%%%%%%%%%%%

Suppose now that $(M,a)$ is a $c$-symplectic manifold, that is, a 
closed  oriented $2n$-dimensional
manifold equipped with a class $a\in H^2(M;\R)$ such that
$a^n > 0$.  Then the analogue of the symplectic group
$\Symp$ is the group $\Diff^a$
 of all diffeomorphisms whose action on 
$H^2(M;\R)$ preserves $a$.  Thus its identity component is the full group $\Diff_0$.  Since this is a simple group, there is no {\it subgroup} corresponding to the Hamiltonian group.  However, as noted in K\c edra--McDuff~\cite{KMc}
there is a covering group, which may be defined as follows.

Consider the
$a$-Flux homomorphism
$
\Flux^a: \pi_1(\Diff_0)\to H^1(M;\R),
$
whose value on a loop $\la$ in 
$\Diff_0(M)$ is  the cohomology class
$\Flux^a(\la)$, where, for a $1$-cycle $\ga$ in $M$,
$\Flux^a(\la)(\ga)$ is given by evaluating 
$a$ on the 
$2$-cycle defined by the map $\T^2\to
 M,(s,t)\mapsto \la(t)\bigl(\ga(s)\bigr)$. 
Set\footnote
{
By McDuff~\cite[Thm~6]{Mc} $\Ga_a$ need not be discrete.
}
$$
\Ga_a: = {\rm Im} (\Flux_a)\subset H^1(M;{\mathcal P}_a^\Z).
$$
We define the $c$-Hamiltonian group
${\rm HDiff}^a_0$ to be the corresponding
covering space of $\Diff_0: = \Diff_0(M)$.
Thus there is a group extension: 
$$
\Ga_a\to {\rm HDiff}_0^a\to\Diff_0.
$$

A similar construction in the symplectic case gives
 a group extension 
$$
\Ga \to \Hh_0\to \Symp_0,
$$
such that the inclusion $\Ham\hookrightarrow \Symp_0$ has a natural 
lift to
a homotopy equivalence  $\io:\Ham\to \Hh_0$.  To see this, note that
if $h\in \Ham$, the element in $\Hh_0$ represented by
the pair $(h, \ga)$, where
 $\ga$ is any path in $\Hh$ 
from the identity to $h$, is independent of the choice of  $\ga$ because
 the difference between two such paths 
lies in the 
kernel of $\Flux^{[\om]}$.

Here are some analogs for $\Diff^a$ of the
  questions considered earlier.

\begin{question}\labell{qu1} Suppose that $M\to P\to B$ 
is a bundle with structural group
$\Diff^a$.  When does the class $a\in H^2(M)$ extend to a class 
$\Tilde a\in H^2(P)?$
\end{question}
 
\begin{question}\labell{qu2} When is there a group extension 
$$
 \Ga_a\to {\rm HDiff}^a\to \Diff^a
 $$
 that restricts over $\Diff_0$ to the extension
 $
 \Ga_a\to  {\rm HDiff}^a_0 \to \Diff_0?
 $
 \end{question}
 
Neither question has an obvious answer.  It is also not clear what
relation they have to each other.  For example, consider the case 
when $\Ga_a: = \im (\Flux^a) = 0$.
 Then the cover is trivial and so always extends, but  $a$
may not always extend. We shall see below that Question~\ref{qu2}
is, at least
to some extent, analogous to Question~\ref{qu:2} which asks when
Flux extends.  However, in the present situation there is so far no analog of Theorem~\ref{thm:obs} which shows the close relation between the existence of
an extension of Flux and the obstruction cocycle $\OM$.

%We now show that
%these questions do have some relation to our previous discussion.

The cocycle $\eps_0\in  H^2_{\EM}(\Diff_0^a; \Ga_a)$
that determines the extension  $
 \Ga_a\to  {\rm HDiff}^a_0 \to \Diff_0
 $  is defined as follows.
 For each $g\in \Diff_0$, choose a path $\ga_{g}$ from the identity 
 element to $g$ and then define $\eps_0(g,h)$ to be the value of 
Flux$^a$
 on the loop formed by going along $\ga_h$ to $h$ then along $\ga_gh$ 
 to $gh$ and then back along $\ga_{gh}$. 
In the symplectic case $\eps_0(g,h)$ can be defined as the sum: 
 $$
 \eps_0(g,h) = \Flux(\ga_h) + \Flux(\ga_g) - \Flux(\ga_{gh}).
 $$
The next result shows that $\eps_0$ extends to a cocycle $\eps$ on $\Symp$
when Flux extends. 

%Hence Question~\ref{qu2}, which asks when  $\eps_0$ 
%extends to $\eps\in H^2_{\EM}(\Diff; \Ga_a)$, is, at least
%to some extent,
%the analog of  
%Question~\ref{qu:2}. However, in the present situation there is so far no 
%analog of Theorem~\ref{thm} which shows the close relation between the 
%existence of
%an extension of Flux and the obstruction cocycle $\OM$.

 \begin{prop}\labell{le:eps} Suppose that $\Flux$ extends to a crossed
 homomorphism $\TF:\Symp\to H_\R/\Ga$, and set $\Hh_{\Tilde F}: = \ker
 \Tilde F$.  For each $g\in \Symp$, pick an element $x_g\in \Hh_{\Tilde
 F} $ that is isotopic to $g$ and choose a path $\ga_{g}$ in $\Symp_0$
 from the identity to $gx_g^{-1}$.  Then: \sss
     
\NI {\rm (i)}\,The formula
 \begin{equation}\labell{eq:eps}
  \eps(g,h) = \Flux(\ga_h) + h^*\Flux(\ga_g) - \Flux(\ga_{gh}),
\end{equation}
defines an element
     $[\eps]\in H^2_{\EM}(\Symp;\Ga)$ that is independent of choices.
\sss
 
 \NI
 {\rm (ii)}   The inclusion $\Hh_{\TF}\to \Symp$ lifts to a homotopy 
 equivalence between $\Hh_{\TF}$ and the covering group $\Tilde\Gg$ of $\Symp$ 
defined by $\eps$.
 \end{prop}
\begin{proof}  To check (i) first observe that 
     $x_g$ is determined modulo an element in $\Ham$ and so
     the element 
$\Flux(\ga_h)$ is independent of choices modulo an element $c_h\in 
\Ga$ that depends on the choice of path $\ga_h$.
The three paths $\ga_h, \ga_g hx_h^{-1}, -\ga_{gh}$ no longer make up 
a 
triangle; to close them up into a loop one needs to add a path from 
$gx_g^{-1}hx_h^{-1}$ to $gh x_{gh}^{-1}$.  But there are elements 
$k,k'\in \Ham$ such that
$$
gx_g^{-1}hx_h^{-1} = gh x_g^{-1}k x_h^{-1} = gh x_{gh}^{-1} k'
$$
and so we can choose this path to lie in $\Ham$. It follows that 
formula~(\ref{eq:eps}) does define an element $\eps(g,h)\in \Ga$.
The cocycle condition
$$
\de \eps(g,h,k) =\eps (h,k) - \eps (g,hk) + \eps(gh,k) - k^*\eps(g,h)
$$
follows by an easy calculation.  Moreover
different choices of the paths $\ga_h$ 
change $\eps$ by a coboundary.
 This proves (i).

To prove (ii) observe that the extension $\Tilde\Gg$ defined by $\eps$ 
has elements  $(g,a)\in \Symp\times \Ga$ where
$$
(g,a) (h,b) = \bigl(gh, h^*a + b + \eps(g,h)\bigr).
$$
Now consider the map $\Tf: \Tilde\Gg\to H_\R$ given by 
$$
\Tf(g,a) = \Tilde\Flux(\ga_g) - a \in H^1(M;\R),
$$
where $\Tilde\Flux(\ga_g)$ is the flux along the path $\ga_g$ in the 
universal cover $\Tilde\Symp$.
It follows immediately from the definitions  that this is a crossed homomorphism with kernel $\Tilde\Hh$, say.  Thus
there is a commutative diagram with exact rows and columns
$$
\begin{array}{ccccc}
&&\Ga &\stackrel{\rm id}\to & \Ga\\
&&\downarrow &&\downarrow\\
\Tilde\Hh &\to &   \Tilde\Gg &\stackrel{\Tf}\to & H_\R\\
\downarrow &&  \downarrow &&\downarrow\\  
\Hh_\TF&\to&  \Symp &\stackrel{\TF}\to & H_\R/\Ga.\end{array}
   $$
Because $H_\R$ is contractible the $5$-lemma implies that the map $\Tilde\Gg\simeq \Tilde\Hh\to \Hh_\Ff$ is
a homotopy equivalence.
\end{proof}

\begin{rmk}\labell{rmk:eps}\rm (i)  
Suppose that $\Poz = \Z$, and repeat the above
argument using the crossed homomorphism $\wF_s^{\Z}$, where $s$ is
$\tau$-compatible.  By Proposition~\ref{prop:Z}, $\Ham^{s\Z}$ maps
onto $\pi_0(\Symp_{\tau})$.  Hence $\Ham^{s\Z}$ is homotopy equivalent
to the covering group of $\Symp_{\tau}$ defined by the cocycle
$\eps_{s}\in H^2_{c\EM}(\Symp_{\tau}; H_\Z)$ of (\ref{eq:eps}).  It follows
from Proposition~\ref{prop:Z} that this 
covering group may be identified with
$\ov\Gg_{\tau}$.
\sss

\NI
(ii)  If $F:G\to \Aa/\La$ is a continuous crossed homomorphism then there is an associated (possibly discontinuous) extension cocycle $\eps_F\in H^2_{EM}(G;\La))$ defined by
$\eps_F(g,h) = h^*f(g) +f(h) - f(gh)$,
where $f:G\to \Aa$ is any lift of $F$.  It is not hard to check that the corresponding extension of $G$ coincides with the covering group $\TG$ constructed from $F$ in Remark~\ref{rmk:uniq}(iii).
\end{rmk}
    
    Observe finally that we can relax Question~\ref{qu2} 
    by asking for an extension of
$\Diff^a$ by the group $H^1(M;{\mathcal P}_a^\Z)$ rather than by
$\Ga$.  If $a$ is a primitive integral class (i.e. ${\mathcal P}_a^\Z
= \Z$), there are natural candidates for such an extension just as in
the symplectic case.  Moreover the existence of these groups
 should have some bearing on Question~\ref{qu1}.
 Since these questions are very similar to those already discussed in
connection with the group $\Ham^{s\Z}$, we shall not pursue them
further here.

%%%%%%%%%%%%%%%%%%%%%%%%%%%%%%%%%%%%%%%%%%%%%%%%%%%%%%%%%%%%%%%%%%%%%%%%

\end{document}